\documentclass[12pt,a4paper]{amsart}
\usepackage{amscd}
\usepackage{amssymb}
\usepackage[centering,text={14.5cm,22cm}]{geometry}
\usepackage{graphicx}
\usepackage{color}
\usepackage[all]{xy}
\usepackage{mathrsfs}
\usepackage{marvosym}
\usepackage{stmaryrd}
\usepackage{srcltx}
\usepackage{hyperref}

\definecolor{shadecolor}{rgb}{1,0.9,0.7}

\newtheorem{theorem}[equation]{Theorem}

\newtheorem{lemma}[equation]{Lemma}

\newtheorem{claim}[equation]{Claim}
\newtheorem{proposition}[equation]{Proposition}
\newtheorem{corollary}[equation]{Corollary}

\newtheorem{definition-lemma}[equation]{Definition-Lemma}

\theoremstyle{definition}
\newtheorem{definition}[equation]{Definition}
\newtheorem{question}[equation]{Question}

\newtheorem{example}[equation]{Example}

\theoremstyle{remark}
\newtheorem{remark}[equation]{Remark}

\numberwithin{equation}{section}
\newtheorem{say}[equation]{}

\newcommand{\EE} {\mathbb{E}}
\newcommand{\HH} {\mathbb{H}}

\newcommand{\ZZ} {\mathbb{Z}}
\newcommand{\QQ} {\mathbb{Q}}

\newcommand{\CC} {\mathbb{C}}

\newcommand{\PP} {\mathbb{P}}
\renewcommand{\AA} {\mathbb{A}}

\newcommand{\coim} {\operatorname{coim}}

\newcommand {\shC}  {\mathcal{C}}
\newcommand {\shD}  {\mathcal{D}}

\newcommand {\E}  {\mathcal{E}}
\newcommand {\shF}  {\mathcal{F}}

\newcommand {\shH}  {\mathcal{H}}

\newcommand {\shI}  {\mathcal{I}}

\newcommand {\shL}  {\mathcal{L}}
\renewcommand {\L} {\shL}

\renewcommand {\SS}  {\mathcal{S}}

\newcommand {\shX}  {\mathcal{X}}
\newcommand {\shY}  {\mathcal{Y}}

\newcommand {\m}  {\mathfrak{m}}

\newcommand {\Char} {\operatorname{char}}

\newcommand {\coker} {\operatorname{coker}}

\newcommand {\Def} {\operatorname{Def}}

\newcommand {\dual} {\vee}

\newcommand {\Ext}  {\operatorname{Ext}}
\newcommand {\ext}  {\operatorname{Ext}}
\newcommand {\lext}  {\underline{\operatorname{Ext}}}

\newcommand {\Hom}  {\operatorname{Hom}}
\renewcommand {\hom}  {\operatorname{Hom}}
\newcommand {\lhom} {\underline{\operatorname{Hom}}}

\newcommand {\hyperext} {\EE\!\text{\textit{xt}}}

\newcommand {\im}  {\operatorname{im}}

\renewcommand {\ker } {\operatorname{ker}}

\newcommand {\loc} {\mathrm{loc}}

\renewcommand{\O}  {\mathcal{O}}

\renewcommand{\P}  {\PP}

\newcommand {\Pic}  {\operatorname{Pic}}

\newcommand {\pr}  {\operatorname{pr}}

\newcommand {\scrB}  {\mathscr{B}}
\newcommand {\B} {\mathscr{B}}
\newcommand {\scrC}  {\mathscr{C}}

\newcommand {\Supp} {\operatorname{Supp}}

\newcommand {\Sing} {\operatorname{Sing}}

\newcommand {\Spec} {\operatorname{Spec}}

\newcommand {\Tor}  {\operatorname{Tor}}

\newcommand {\D} {\mathscr{D}}
\newcommand {\T} {\Theta}
\newcommand {\X} {\shX}

\newcommand {\I} {\shI}

\newcommand {\Pone} {\PP^1}
\newcommand {\Pthree} {\PP^3}
\newcommand {\Pfour} {\PP^4}
\newcommand {\Ptwo} {\PP^2}

\def\dual#1{{#1}^{\scriptscriptstyle \vee}}
\def\mapright#1{\smash{
  \mathop{\longrightarrow}\limits^{#1}}}
\def\mapleft#1{\smash{
  \mathop{\longleftarrow}\limits^{#1}}}
\def\exact#1#2#3{0\to#1\to#2\to#3\to0}

\def\mydate{\ifcase\month \or January\or February\or March\or
April\or May\or June\or July\or August\or September\or October\or 
November\or December\fi \space\number\day,\space\number\year}

\begin{document}

\title{Deforming Calabi-Yau threefolds}

\author{Mark Gross}

\date{\today}
\maketitle

\setcounter{tocdepth}{2}
\tableofcontents

%=================================================================================

\section*{Introduction.}\footnote{This is a 2025 update of this paper, 
published
in 1997. This revision addresses two issues. First, it corrects an error
concerning the definition of $T^1$ in the case when the deformation functor
is not pro-representable. This was flagged in \cite{KawaErr}. In particular,
the condition $(H_5)$ given in the published version of the paper is not 
necessary and the relevant functors in any case do not satisfy the condition
$(H_5)$. Second, it provides some details requested by R.\ Friedman for
compatibility of certain exact sequences in the proof of Theorem~\ref{2.2}.
The results of the paper remain unchanged.}

There are two things that quickly become clear in surveying the work done on the
subject of Calabi-Yau threefolds, by which we mean a projective threefold $X$ (with some
specified class of possible singularities) with $K_X=0$ and $h^1(\O_X)=0$. 
First, there are a huge number of such threefolds, even non-singular. Second, one
reason there appears to be so many is this: suppose you degenerate a non-singular
Calabi-Yau
$X$ to a threefold $X'$ with canonical singularities. If $X'$ has
a crepant desingularization $\tilde X$, then  $\tilde X$ will not, in general, be
in the same deformation family as $X$ or even be diffeomorphic to 
$X$. As a result, an
innocent-enough Calabi-Yau threefold such as the quintic hypersurface in 
$\Pfour$
can have hundreds of degenerations, if not more, with crepant resolutions, thus
giving rise to huge numbers of other Calabi-Yaus. 

Classification of these degenerations seems to be a hopeless problem. All we
know, of course, is that there are only a finite number of families of
such degenerate quintics. So, if we want to try to simplify the
task of Calabi-Yau classification, it might help to concentrate
on Calabi-Yau threefolds which do not arise as crepant resolutions of degenerations
of other Calabi-Yau threefolds.  This motivates the following definition:

\begin{definition}
A non-singular Calabi-Yau threefold $\tilde X$ is 
{\it primitive} if there is no birational contraction $\tilde
X\rightarrow X$ with $X$ smoothable to a Calabi-Yau threefold which is not
deformation equivalent to $\tilde X$. 
\end{definition}

(The last condition is included to rule out the sort of possibility that can occur,
say, if one has contracted an elliptic scroll to a curve. See \cite[Ex.~4.6]{[50]} for such an example.)

While I will not study primitive Calabi-Yau threefolds directly in this
paper (I defer this until [12]), this definition will still guide our inquiry.
In particular, in order to understand which Calabi-Yaus are not
primitive, we should understand the answer to the following question: given a
Calabi-Yau threefold
$\tilde X$, when can we find a birational contraction morphism $\tilde X\rightarrow
X$ such that $X$ is smoothable? Now, any birational contraction will yield a
threefold
$X$ with at worst canonical singularities. So if we want to begin to understand
the smoothability of such threefolds, we first need to understand if they
have obstructed deformation theory or not. Thus we have

\begin{question}
Given a Calabi-Yau threefold $X$ with canonical
singularities, is $\Def(X)$ non-singular? If not, can we get some
reasonable dimension estimates for components of $\Def(X)$?
\end{question}

As already shown in [11], if $X$ has canonical singularities, $\Def(X)$ can indeed be
singular. However, as we shall show in this paper, we can still control 
the dimension of components of $\Def(X)$ if
$X$ has canonical singularities. The principle we discover is that obstructions
to deforming $X$ are essentially the obstructions to deforming a germ of the
singularities of $X$. This gives a further generalization of the
Bogomolov-Tian-Todorov unobstructedness theorem. This material is covered in
\S 2, with preliminaries in \S 1. We generalise the techniques of \cite{[28]} 
and
some ideas from deformation theory due to Z. Ran.

Now, in the attempt to understand when Calabi-Yau threefolds with canonical
singularities are smoothable, it will certainly be hopeless to try to understand all
possible canonical singularities which can arise and then determine which ones are
smoothable. Nevertheless, we still  obtain strong results. We have

\begin{theorem} [Theorem \ref{3.8}] 
Let $\tilde X$ be a non-singular Calabi-Yau threefold,
and $\pi:\tilde X\rightarrow X$ be a birational contraction morphism,
such that $X$ has isolated complete intersection singularities. Then
there is a deformation of $X$ which smooths all singular points of $X$
except possibly the ordinary double points of $X$.
\end{theorem}

Results of Namikawa and Steenbrink \cite{[30]} then allow us to extend this result to the case
that
$\tilde X$ has terminal singularities.
 
In the case that the isolated singularity is not complete intersection, we
have a much weaker statement. The hypothesis that $X$ is $\QQ$-factorial
is necessary to ensure that $X$ has enough infinitesimal deformations. Furthermore, 
since at the moment we do not have any real control over how bad the deformation
space of a canonical singularity can be, we include a rather artificial hypothesis on
the singularities we will consider (see Definition \ref{4.2}) called {\it good}.
Hopefully some of these hypotheses can be removed at a future date. We show in \S 5
that some simple classes of singularities are good, and this is enough for initial
applications of our results. We have

\begin{theorem}
[Theorem \ref{4.3}] Let $\tilde X$ be a non-singular Calabi-Yau threefold and
$\pi:\tilde X\rightarrow X$ a birational contraction such that $X$ is 
$\QQ$-factorial and for each $P\in Sing(X)$, the germ $(X,P)$ is good.
Then $X$ is smoothable.
\end{theorem}

Recall that a birational projective contraction $\pi:\tilde X\rightarrow X$
is {\it primitive} if it cannot be factored in the projective category. One
application of Theorem \ref{4.3} given in \S\ref{sec5} is

\begin{theorem} [Theorem \ref{5.8}]
Let $\pi:\tilde X\rightarrow X$ be a primitive contraction
contracting a divisor $E$ to a point. Then $X$ is smoothable unless $E\cong\Ptwo$
or $F_1$.
\end{theorem}

Finally, a bit of history on these questions. The unobstructedness question has been
answered positively for varying degrees of singularities: for non-singular
Calabi-Yaus of any dimension by Bogomolov, Tian \cite{[45]} and Todorov 
\cite{[47]}, with 
algebraic proofs given by Ran \cite{[33]}, Kawamata \cite{[18]}, and Deligne; for Calabi-Yaus
with ordinary double points by Kawamata \cite{[18]} and Tian \cite{[46]}; for Calabi-Yaus with
Kleinian singularities and orbikleinfold singularities by Ran \cite{[34],[37]}; and finally
for Calabi-Yau threefolds with rational isolated complete intersection singularities
by Namikawa \cite{[28]}. Results on smoothability of singular Calabi-Yau threefolds were
first obtained by Friedman \cite{[8]} for Calabi-Yau threefolds with ordinary double
points, and by \cite{[30]} for Calabi-Yaus with hypersurface singularities. There is some
overlap between Theorem \ref{3.8} for hypersurface singularities and the results of \cite{[30]},
using similar methods.  I received a preprint version of \cite{[30]} just as I was finishing
the first version of this paper.

I
would like to thank H.\ D'Souza, J.\ Koll\'ar, Z.\ Ran, M.\ Stillman, D.\ Van Straten,
P.M.H.\ Wilson and B.\ W\"ubben for useful discussions, suggestions and answers to questions during
this work. I thank R.\ Friedman for the interest and discussion leading to
the current revision.

\section{Some Deformation Theory}
\label{sec1}

There is little that is new in this section, but I will be needing a number of
variants on results primarily due to Ran and Kawamata. I give complete
proofs here of the precise statements I will  be needing.
I begin by reviewing some facts of deformation theory.

\begin{say}
\label{(1.1)}
The general context from \cite{[41]} is as follows. Let $k$
be a field and let $\Lambda$ be
a complete Noetherian local $k$-algebra with residue field $k$ and maximal ideal
${\bf m}_{\Lambda}$. We denote by
${\bf C}_{\Lambda}$ the category of Artin local $\Lambda$-algebras with residue field
$k$ with local homomorphisms. We are interested in deformation
functors $D:{\bf C}_{\Lambda}\rightarrow {\bf Ens}$ (here ${\bf Ens}$ is the category
of sets), with $D(k)$ consisting of one element. 
\end{say}

\begin{say}
\label{(1.2)}
We will only consider functors $D$ which are pro-representable
or have a hull (miniversal space). Thus there is a complete
local
$\Lambda$-algebra $S$ with residue field $k$ and a morphism of functors
$\hom(S,\cdot)\rightarrow D$. This morphism is an isomorphism if $D$ is
pro-represented by $S$. If $S$ is a hull of $D$, this morphism is only smooth and
induces an isomorphism on tangent spaces. Here
$\hom$ denotes local
$\Lambda$-algebra homomorphisms.
We will always write $S\cong R/J$, where $R=\Lambda[[x_1,\ldots,x_r]]$
with maximal ideal ${\bf m}_{R}={\bf m}_{\Lambda}R+
(x_1,\ldots,x_r)$, and $J\subseteq {\bf m}_{\Lambda}R+
{\bf m}_R^2$  an ideal.
\end{say}

\begin{say}
\label{(1.3)}
Following \cite{[18]}, set
\begin{align*}
A_n={} & k[t]/(t^{n+1})\\
B_n={} & A_n\otimes_k A_1=k[x,y]/(x^{n+1},y^2)\\
C_n= {} & B_{n-1}\times_{A_{n-1}} A_n=k[x,y]/(x^{n+1},x^ny,y^2)
\end{align*}
and let $\alpha_n:A_{n+1}\rightarrow A_n$, $\beta_n:B_n\rightarrow A_n$,
$\gamma_n:B_n\rightarrow C_n$, and $\xi_n:B_n\rightarrow B_{n-1}$ be the natural
maps. Define
$\epsilon_n:A_{n+1}\rightarrow B_n$ by $t\mapsto x+y$ and $\epsilon_n':A_n\rightarrow
C_n$ by $t\mapsto x+y$ also.
\end{say}

\begin{say}
\label{(1.4)}
Now consider the case that $\Lambda=k$ and that $D$ is prorepresentable. For
$X_n\in D(A_n)$ we define the first order tangent space  of $X_n$, 
$$T^1(X_n/A_n)=\{Y_n\in D(B_n)\,|\,D(\beta_n)(Y_n)=X_n\}.$$
This set has a natural $A_n$-module structure
as follows. First, if
$\alpha,\beta\in T^1(X_n/A_n)$, we need to define $\alpha+\beta$. We have
$\alpha\times \beta\in D(B_n)\times_{D(A_n)} D(B_n)=D(B_n\times_{A_n} B_n)$ by
pro-representability. We have an isomorphism of $k$-algebras
\[
B_n\times_{A_n} B_n\cong
k[x,y,y']/(x^{n+1},y^2,yy',y'^2),
\] 
and there is a natural map 
$B_n\times_{A_n} B_n\rightarrow B_n$ via $x\mapsto x, y\mapsto y, y'\mapsto y$. 
Then
$\alpha+\beta$ is the image of $\alpha\times\beta$ in $D(B_n)$ under this map.
It is easy to check that $\alpha+\beta\in T^1(X_n/A_n)$. 

Secondly, if $a\in A_n$, we have the endomorphism $a:B_n\rightarrow B_n$ given
by $x\mapsto x$, $y\mapsto ay$. Then if $\alpha\in T^1(X_n/A_n)$, so is
$a\alpha=D(a)(\alpha)$.
\end{say}

\begin{say}
\label{(1.5)} 
In this paper, we will on occasion work with functors which are not necessarily
pro-representable, but which do have a hull. The deformation functor $D$ has
a hull if and only if Schlessinger's conditions $(H_1)$-$(H_3)$ of 
\cite[Thm.~2.11]{[41]}
are satisfied. We recall these here. For any morphisms $A'\rightarrow A$,
$A''\rightarrow A$ in ${\bf C}_{\Lambda}$, there is a natural functorial map
$$\psi_{A',A'',A}:D(A'\times_A A'')\rightarrow D(A')\times_{D(A)} D(A'').$$
A surjective map $A'\rightarrow A$ in ${\bf C}_{\Lambda}$ is a {\it small extension}
if its kernel is a non-zero principal ideal $(t)$ such that ${\bf m}_{A'}(t)=0$.
 The conditions $(H_1)$-$(H_3)$ are:
\begin{itemize}
\item[$(H_1)$] $\psi_{A',A'',A}$ is a surjection whenever $A''\rightarrow A$ is a
small extension.
\item[$(H_2)$] $\psi_{A',A'',A}$ is a bijection when $A=k$, $A''=A_1$. 
\item[$(H_3)$] $\dim_k(T^1(X/k))<\infty$.
\end{itemize}

$(H_2)$ guarantees that $T^1(X/k)$ has a $k$-vector space structure as in \ref{(1.4)}.
If in addition
$$\hbox{$(H_4)$ $\psi_{A',A',A}$ is a bijection for any small extension
$A'\rightarrow A$\hfill}$$
is satisfied then $D$ is pro-representable. (This is also a part of 
\cite[Thm.~2.11]{[41]}).
\end{say}

\begin{say}
As pointed out in \cite{KawaErr}, the definition for $T^1(X_n/A_n)$
given in \ref{(1.4)} is not correct if $D$ is not pro-representable,
and in particular need not be a vector space.\footnote{The original version
of this paper tried to get around this issue by introducing a property
$(H_5)$. However, in (1.6) of the published version of this paper,
the functor $D_X$ does not satisfy condition $(H_5)$. The construction
given is not well-defined since the result depends on choices.
Indeed, if $X_{A'}$, $X_{A''}$, $X_A$ are deformations
of $X$, $X_{A'}\otimes_{A'} A\mapright{\phi_1} X_A\mapleft{\phi_2} X_{A''}
\otimes_{A''} A$, with $\phi_1$, $\phi_2$ isomorphisms, then
the isomorphism class of
$(X,\O_{X_{A'}}\times_{\O_{X_A}} \O_{X_{A''}})$ will depend on the particular choice
of the isomorphisms $\phi_1$, $\phi_2$,
as an automorphism of $X_A$ may not lift
to an automorphism of $X_{A'}$ or $X_{A''}$.}
Instead, following \cite{KawaErr}, one may still define an appropriate
space if $D$ is defined as the functor of flat deformations of a ringed
space $(X,\O_X)$ defined over $\Spec k$. However, in this revision
we will take a more general view-point and instead use the language
of categories cofibred in groupoids as developed in \cite{Rim} and
exposited in \cite[Tag~06G7]{Stacks}. We refer the reader to 
the latter reference for the principle definitions. 

A category cofibred in groupoids is a category $\shD$ equipped with a functor
$p:\shD\rightarrow \shC_{\Lambda}$ such that $p^{\mathrm{opp}}:
\shD^{\mathrm{opp}}\rightarrow \shC_{\Lambda}^{\mathrm{opp}}$ 
is a category fibred
in groupoids. We write $\shD(A)$ for the sub-category of
$\shD$ over $A\in \shC_{\Lambda}$; this is a groupoid. Further,
the definition of category fibred in groupoids tells us that
for a morphism $f:A\rightarrow B$ in $\shC_{\Lambda}$, there is for
each $X$ an object in $\shD(A)$ a push-forward object $f_*X$ in 
$\shD(B)$ equipped with a morphism $X\rightarrow f_*X$ over $f$.
This push-forward object is unique up to unique isomorphism 
(see \cite[Tag 06SH]{Stacks}). Further, push-forward has the property
\begin{equation}
\label{eq:pushfoward property}
\substack{\hbox{Given $f:A\rightarrow B$, $g:B\rightarrow C$, there is
a unique morphism $f_*X \rightarrow (g\circ f)_*X$,}\\
\hbox{ and hence a unique
isomorphism $g_*f_*X \rightarrow (g\circ f)_* X$.}}
\end{equation}
This again follows
from the properties of categories cofibred in groupoids, see 
\cite[Tag 06GJ]{Stacks}.

We use the
convention that we write $D(A)$ for the set of isomorphism classes
of $\shD(A)$. This gives a functor $D:\shC_{\Lambda}\rightarrow
\mathbf{Ens}$ given by $A\mapsto D(A)$.

\begin{example}
\label{ex:principal}
The principal example for us is given by deformations of a ringed space
$(X,\O_X)$ defined over $\Spec k$. Objects are data $(X_A,\O_{X_A})$
flat over $\Spec A$ and an isomorphism $\phi_{X_A}:X_A\otimes_A k\cong
X$. (Here we write $(Y,\O_Y) \otimes_A B$ for $(Y,\O_Y\otimes_A B)$.)
A morphism $X_A\rightarrow X_B$ lying over a homomorphism $A\rightarrow B$
is an isomorphism $X_A\otimes_A B \cong X_B$ compatible with the
isomorphisms $\phi_{X_A},\phi_{X_B}$. This is an example of a predeformation
category, see \cite[Tag 06GS]{Stacks}, in that $\shD(k)$ is equivalent
to a category with a single object and a single morphism, i.e., the identity.
So in particular, the corresponding functor $D:\shC_{\Lambda}\rightarrow
{\bf Ens}$ is a deformation functor.
\end{example}

The analogue of the Schlessinger conditions $(H_1)$,$(H_2)$ is the
Rim-Schlesssinger condition (RS), see \cite[Tag 06J3]{Stacks}. 
We say $\shD$ satisfies (RS) if for ever diagram in $\shD$
\begin{equation}
\label{eq:empty square}
\xymatrix@C=20pt
{
&X_{A_2}\ar[d]\\
X_{A_1}\ar[r]&X_A
}
\end{equation}
lying over
\[
\xymatrix@C=20pt
{
&A_2\ar[d]\\
A_1\ar[r]&A
}
\]
in $\shC_{\Lambda}$ with $A_2\rightarrow A$ surjective, there exists
a fiber product $X_{A_1}\times_{X_A} X_{A_2}$ in $\shD$ such that the
diagram
\[
\xymatrix@C=20pt
{
X_{A_1}\times_{X_A} X_{A_2}\ar[r]\ar[d]&X_{A_2}\ar[d]\\
X_{A_1}\ar[r]&X_A
}
\]
lies over
\begin{equation}
\label{eq:full square}
\xymatrix@C=20pt
{
A_1\times_A A_2\ar[r]\ar[d]&A_2\ar[d]\\
A_1\ar[r]&A
}
\end{equation}

\begin{remark}
\label{rem:06J3}
We note from \cite[Tag~06J3]{Stacks} that if $\shD$ satisfies (RS)
and $A_2\rightarrow A$ is surjective, then any commutative square
completing \eqref{eq:empty square} and lying over \eqref{eq:full square}
is a fibre square.
\end{remark}

\begin{example}
Example~\ref{ex:principal} satisfies condition (RS), with
\[
(X_{A_1},\O_{X_{A_1}})\times_{(X_A,\O_{X_A})} (X_{A_2},\O_{X_{A_2}})
=(X, \O_{X_{A_1}}\times_{\O_{X_A}} \O_{X_{A_2}}).
\]
\end{example}

%For $A\in {\bf C}_k$, an infinitesimal deformation of $X$ over $A$ is a ringed space
%$(X,\O_{X_A})$ such that $\O_{X_A}$ is a sheaf of flat $A$-algebras and
%$\O_{X_A}\otimes_A k=\O_X$. $(X,\O_{X_A})$ and $(X,\O_{X_A}')$ are isomorphic
%deformations
%if there is an isomorphism $\O_{X_A}\rightarrow \O_{X_A}'$ which is the identity
%upon tensoring with $k$. We define $D_X:{\bf C}_k\rightarrow {\bf Ens}$ by
%$$D_X(A)=\{\hbox{Isomorphism classes of deformations of $X$ over $A$.}\}.$$
Following \cite{KawaErr}, one may define $T^1(X_n/A_n)$ in this more
general context. Let $\shD$ be a predeformation category satisfying
(RS) and $X_n$
an object of $\shD(A_n)$. Then $T^1(X_n/A_n)$ is defined to be
the set of isomorphism classes of pairs $(Y_n,\psi_n)$ consisting
an object $Y_n\in \shD(B_n)$ and a morphism
$\psi_n:Y_n \rightarrow X_n$ over $\beta_n:B_n\rightarrow A_n$.  

Similarly to \ref{(1.5)}, $T^1(X_n/A_n)$
has a natural $A_n$-module structure. Indeed, given $(Y_n,\psi_n),
(Y_n',\psi_n')\in T^1(X_n/A_n)$, $Y_n\times_{X_n} Y_n'$ is well-defined
up to isomorphism in $\shD(B_n\times_{A_n} B_n)$. Here the two maps
$Y_n,Y_n'\rightarrow X_n$ used to define the fibre product are $\psi_n$ and
$\psi'_n$ respectively.

The natural map $s:B_n\times_{A_n} B_n
\rightarrow B_n$ as before then yields an object $Y_n''=
s_*(Y_n\times_{X_n} Y_n')$ of $\shD(B_n)$.
Further, $Y_n\times_{X_n} Y_n'$ carries a well-defined morphism to $X_n$,
so $Y_n''$ does also by the cofibred category property, see condition
(2) of \cite[Tag 06GJ]{Stacks}.
Multiplication by elements of $A_n$, is accomplished in the
same way as in \ref{(1.4)}.
\end{say}

\begin{say}
\label{(1.7)} 
If $D$ is a deformation functor on ${\bf C}_{\Lambda}$, we say that a
$k$-vector space $T^2$ is an obstruction space for $D$ if
whenever we have a surjection $\phi:A'\rightarrow A$
in ${\bf C}_{\Lambda}$ with $I=\ker\phi$ annihilated by the maximal ideal of $A'$, 
we get a sequence
\[
\xymatrix@C=20pt
{
D(A')\ar[r]^{D(\phi)}& D(A) \ar[r]^{\delta}& T^2\otimes I,
}
\]
and this sequence is exact in the sense that if
$\alpha\in D(A)$, then 
$\delta(\alpha)=0$ if and only if $\alpha$ is in the image
of $D(\phi)$. Furthermore, the obstruction map should be functorial, so given in
addition $\phi':B'\rightarrow B$ surjective, $I'=\ker\phi'$ annihilated by the
maximal ideal of $B'$ and a commutative diagram
\[
\xymatrix@C=20pt
{
A'\ar[r]^{\phi}\ar[d]_{\beta'}&A\ar[d]^{\beta}\\
B'\ar[r]_{\phi'}&B
}
\]
there is a commutative diagram
\[
\xymatrix@C=20pt
{D(A')\ar[r]^{D(\phi)}\ar[d]_{D(\beta')}&D(A)\ar[r]^{\delta}\ar[d]_{D(\beta)}
&T^2\otimes I\ar[d]^{1_{T^2}\otimes\beta'}\\
D(B')\ar[r]_{D(\phi')}&D(B)\ar[r]_{\delta'}
&T^2\otimes I'}
\]

If $D$ is prorepresentable by $S\cong R/J$ as in \ref{(1.2)}, it is easy to describe the
obstruction theory of
$D$. Let $T^2$ be the $k$-vector space $(J/{\bf m}_RJ)^{\vee}$. 
The obstruction map $\delta$ can be described as follows.
Given $f\in D(A)=\hom(R/J,A)$, let $f(x_i)=\alpha_i\in A$. Choose any
lifting of $\alpha_i$ to $\alpha_i'\in A'$; this defines a $\Lambda$-algebra
homomorphism
$f':R\rightarrow A'$. Now given an element $\beta\in J$, we must have $f'(\beta)\in
I$, since
$f(\beta)=0$. Furthermore, if $\beta\in {\bf m}_RJ$, then $f'(\beta)=0$, since
$f'$ is a local homomorphism and $I$ is annihilated by the maximal ideal of $A'$. 
Thus
$f'$ induces a
$k$-vector space map
$J/{\bf m}_RJ\rightarrow I$, i.e. an element of $T^2\otimes I$. We then define
$\delta(f)$ to be this element. It is easy to see that $\delta(f)$ does not depend
on the choice of the lifting. Furthermore, it is clear that $f'$ induces a map
$f':R/J\rightarrow A'$ if and only if $\delta(f)=0$. Note that if one choice of the $\alpha_i'$ 
provides a lifting $f':R/J\rightarrow A'$, then
any choice does, and the set of possible liftings is a principal homogeneous space
over the vector space $(\m_R/(\m_{\Lambda}R+\m_R^2))^{\vee}\otimes I$. 
 
We note also that $(J/{\bf m}_RJ)^{\vee}$ is naturally isomorphic to $T^1(S/\Lambda,
k)$ by \cite[3.1.2]{[25]}, where here $T^1$ is the first cotangent functor of Lichtenbaum
and Schlessinger.
\end{say}

The following is the version of Ran's and Kawamata's statements of
Ran's $T^1$-lifting criterion (\cite{[35]} and \cite{[18]}) as
advertised in \cite{KawaErr}, enhanced to the setting of 
predeformation categories. 

\begin{theorem}[Ran,Kawamata]
\label{1.8}
Let $k$ be a field of characteristic 0, $\shD\rightarrow\shC_{k}$
a predeformation category satisfying (RS), 
$D$ the the corresponding deformation functor. Suppose that $D$
has an obstruction space $T^2$. Then for each
$X_n\in \shD(A_n)$, $X_{n-1}=\alpha_{n-1,*}X_n$, let 
$\alpha\in T^1(X_{n-1}/A_{n-1})$
be given by the pair $(\epsilon_{n-1,*}X_n, \psi_{n-1})$ where 
$\psi_{n-1}$ is the canonical morphism 
$\epsilon_{n-1,*}X_n\rightarrow X_{n-1}$ guaranteed by 
\eqref{eq:pushfoward property}.
Then there exists an
$X_{n+1}\in D(A_{n+1})$ with $D(\alpha_n)(X_{n+1})=X_n$ if and only if $\alpha$ is in the image
of the natural map
$$\xi_{n*}:T^1(X_n/A_n)\rightarrow T^1(X_{n-1}/A_{n-1}).$$
\end{theorem}

\begin{proof}
This proof is simply a very minor modification of Kawamata's proof
of Ran's $T^1$-lifting criterion in \cite{[18]}. 

First, the map $\xi_{n*}$ is defined as follows. An element of
$(Y_n,\psi_n)$ of $T^1(X_n/A_n)$ induces a composed map
$Y_n\rightarrow X_n\rightarrow \alpha_{n-1,*}X_n=X_{n-1}$, and hence
a morphism $\psi_{n-1}:Y_{n-1}:=\xi_{n*}Y_n\rightarrow X_{n-1}$ by
\eqref{eq:pushfoward property}.

As in \cite[pg.~185]{[18]}, we have
a commutative diagram with exact rows: 
\[
\xymatrix@C=30pt
{
D(A_{n+1})\ar[r]^{D(\alpha_n)}
\ar[d]_{D(\epsilon_n)}
&D(A_n)\ar[r]^{\delta_1}
\ar[d]_{D(\epsilon_{n}')}
&T^2\otimes (t^{n+1})
\ar[d]^{1_{T^2}\otimes
\epsilon_n}\\ 
D(B_n)\ar[r]_{D(\gamma_n)}&D(C_n)\ar[r]_{\delta_2}
&T^2\otimes (x^ny)}
\]
Note that $1_{T^2}\otimes\epsilon_n$ is an isomorphism 
(here we need $\Char k=0$, since $\epsilon_n(t^{n+1})=(n+1)x^ny$ in $B_n$).

Put $Y_{n-1}:=\epsilon_{n-1,*}X_n$, so $\alpha=(Y_{n-1},\psi_{n-1})$.
Note that we have a commutative diagram, again using 
\eqref{eq:pushfoward property},
\[
\xymatrix@C=30pt
{
\epsilon_{n*}'X_n\ar[r]\ar[d]& \pr_{2*}\epsilon_{n*}'X_n=X_n\ar[d]\\
Y_{n-1}=\epsilon_{n-1,*}X_n\ar[r]_{\psi_{n-1}}&X_{n-1}=\alpha_{n-1,*}X_n
}
\]
sitting over the fibre square defining $C_n=B_{n-1}\times_{A_{n-1}} A_n$.
Here 
\[
\pr_2:B_{n-1}\times_{A_{n-1}} A_n\rightarrow A_n
\]
denotes the second projection.
Thus by Remark~\ref{rem:06J3}, $\epsilon_{n*}'X_n \cong 
Y_n\times_{X_{n-1}} X_n$.
Now
\begin{itemize}
\item[] $\exists X_{n+1}\in D(A_{n+1})$ such that $D(\alpha_n)(X_{n+1})=X_n$
\item[$\Leftrightarrow$]$\delta_1(X_n)=0$
\item[$\Leftrightarrow$]$\delta_2(D(\epsilon'_n)(X_n))=0$
\item[$\Leftrightarrow$]$\exists Y_n\in D(B_n)$ with
$D(\gamma_n)(Y_n)=D(\epsilon_n')(X_n)$
\item[$\Leftrightarrow$]$\exists$ an object $Y_n$ in $\shD(B_n)$
fitting into a commutative diagram 
\[
\xymatrix@C=30pt
{
\gamma_{n*}Y_n\ar[r]\ar[d] & X_n\ar[d]\\
Y_{n-1}\ar[r]_{\psi_{n-1}}&X_{n-1}
}
\]
sitting over the fibre square defining
$B_{n-1}\times_{A_{n-1}} A_n$ 
\item[$\Leftrightarrow$]$\exists$ $Y_n$ an object of $\shD(B_n)$
with morphisms
$\psi_n:Y_n\rightarrow X_n$ over $\beta_n$ and $\phi:Y_n\rightarrow
Y_{n-1}$ over $\xi_n$ yielding a commutative diagram
\[
\xymatrix@C=30pt
{
Y_n\ar[r]^{\psi_n}\ar[d]_{\phi}&X_n\ar[d]\\
Y_{n-1}\ar[r]_{\psi_{n-1}}&X_{n-1}
}
\]
(Here, we use $\xi_n$ is the composition $\pr_1\circ\gamma_n$ where
\[
\pr_1:C_n=B_{n-1}\times_{A_{n-1}} A_n\rightarrow B_{n-1}
\]
is the first projection, as well as properties of cofibred categories.)
\item[$\Leftrightarrow$]$\exists\alpha'\in T^1(X_n/A_n)$ mapping to $\alpha$.
\end{itemize}
\end{proof}

\begin{definition}
\label{def:predef functor}
Let $p_i:\shD_i\rightarrow \shC_{\Lambda}$ be two predeformation categories.
A functor $F:\shD_1\rightarrow \shD_2$ is \emph{functor of predeformation 
categories} if $p_2\circ F=p_1$ and $F(\alpha_* X)=\alpha_*F(X)$ for
every $\alpha:A\rightarrow B$ in $\shC_{\Lambda}$ and $X$ an object in
$\shD_1(A)$.
\end{definition}

The following theorem is a slightly more specific version of Ran's Theorem 1.1 of
\cite{[36]}.

\begin{theorem}
\label{1.9}
Let $k$ be an algebraically closed\footnote{This hypothesis is
added as in \cite{KawaErr}.} field of characteristic 0, $\shD_1,\shD_2
\rightarrow {\bf C}_k$ two predeformation categories
with a functor of predeformation categories $F:\shD_1\rightarrow
\shD_2$. Suppose $D_1$ is pro-representable by a $k$-algebra $S\cong P/I$, 
$P=k[[x_1,\ldots,x_s]]$ as in \ref{(1.2)}. Suppose $\shD_2$
satisfies condition (RS) and
$D_2$ has a hull $\Lambda$. Let $T^2_1$ and
$T^2_2$ be $k$-vector spaces and $l:T^2_1
\rightarrow T^2_2$ a $k$-vector space map. Denote by
$X$ a representative of the unique isomorphism class
of object of $\shD_1(k)$. Suppose, for all $n$ and for each object
$X_n$ of $\shD_1(A_n)$ inducing an object $X_{n-1}:=\xi_{n*}X_n$ of
$D_1(A_{n-1})$, there is a commutative diagram
\[
\xymatrix@C=30pt
{
T^1_1(X_n/A_n)\ar[r]^{F_*}
\ar[d]_{\xi_{n*}}
&T^1_2(F(X_n)/A_n)
\ar[d]^{\xi_{n*}}\\
T^1_1(X_{n-1}/A_{n-1})\ar[r]_{F_*}
\ar[d]_{\delta_1}
&T^1_2(F(X_{n-1})/A_{n-1})
\ar[d]^{\delta_2}\\
T^2_1\ar[r]_{l}&T^2_2
}
\]
with exact columns and $l|_{\im(\delta_1)}:\im(\delta_1)\rightarrow T^2_2$
injective and where $F_*(Y_n,\psi_n)$ is $(F(Y_n),F(\psi_n))$. Then
$S\cong R/J$, with $R=\Lambda[[x_1,\ldots,x_r]]$ where
\[
r=\dim_k\ker(T_1^1(X/k)\mapright{F_*} T^1_2(F(X)/k)),
\]
and $J\subseteq
\m_{\Lambda}R+\m_R^2$ an ideal. In addition, there is an ideal
$J'\subseteq J$ with
$\Supp(R/J)=\Supp(R/J')$ and ${\bf m}_R/({\bf m}_R^2+J)\cong {\bf m}_R/({\bf
m}_R^2+J')$ such that $J'$ is generated by
$\dim_k\coker (T^1_1(X/k)\mapright{F_*} T^1_2(F(X)/k))$ elements of $R$.
\end{theorem}

\begin{proof}
Since $\Lambda$ is a hull for $D_2$, there is an induced map
$\Lambda\rightarrow S$, unique only up to the induced map $F_*$ on Zariski tangent
spaces. Fix one such map. Now
$$r=\dim_k \coker(\m_{\Lambda}/\m_{\Lambda}^2\mapright{{F_*}^{\vee}}\m_S/\m_S^2),$$
and if we choose elements $\alpha_1,\ldots,\alpha_r\in \m_S$ which along with
$\im\dual{F_*}$ generate $\m_S/\m_S^2$, we can define a map
$R=\Lambda[[x_1,\ldots,x_r]]
\rightarrow S$ by $x_i\mapsto \alpha_i$. This map is surjective, and if its kernel
is $J$, $S\cong R/J$. Furthermore, $J\subseteq\m_{\Lambda}R+\m_R^2$. Let $D_0:
{\bf C}_{\Lambda}\rightarrow {\bf Ens}$ be the functor pro-represented by 
the complete local $\Lambda$-algebra $S=R/J$.
We will first prove:

\begin{claim}
\label{(1.10)} If $V=\coker(T^1_1(X/k)\rightarrow
T^1_2(F(X)/k))$, then $V\subseteq \dual{(J/\m_RJ)}$, and the obstruction map 
$\delta_0$ in
$$D_0(A_{n+1})\mapright{} D_0(A_n)\mapright{\delta_0} \dual{(J/\m_RJ)}\otimes
(t^{n+1})$$ always takes its values in $V\otimes (t^{n+1})$. 
\end{claim}

\begin{proof}
First note the change of rings sequence for the functors
$T^i$ of Lichtenbaum and Schlessinger (\cite[2.3.5]{[25]})
for
$k\rightarrow\Lambda\rightarrow S$
yields the exact sequence
$$T^0(S/k,k)\rightarrow T^0(\Lambda/k,k) \rightarrow T^1(S/\Lambda,k)
\rightarrow T^1(S/k,k),$$
or equivalently
\begin{equation}
\label{(1.11)}
\dual{({\bf m}_S/{\bf m}^2_S)}\rightarrow \dual{({\bf m}_{\Lambda}/{\bf
m}^2_{\Lambda})} \rightarrow
\dual{(J/\m_RJ)}
\mapright{d} \dual{(I/\m_PI)}.\
\end{equation}
Thus $\ker d=V$. Furthermore, the map $d$ is compatible with the obstruction maps.
Specifically, fixing a $\Lambda$-algebra structure on $A_{n+1}$,
hence inducing one on $A_n$, we obtain the diagram
\[
\xymatrix@C=30pt
{D_0(A_{n+1})\ar[r]\ar[d]&D_1(A_{n+1})\ar[d]\\
D_0(A_n)\ar[r]\ar[d]_{\delta_0}&D_1(A_n)\ar[d]_{\delta_1'}\\
\dual{(J/\m_RJ)}\otimes (t^{n+1})\ar[r]_{d\otimes 1}&\dual{(I/\m_PI)}\otimes
(t^{n+1})
}
\]
where the top two horizontal maps are the forgetful inclusions
$\hom_{\Lambda-alg}(R/J,A_i)\subseteq
\hom_{k-alg}(R/J,A_i)$ for $i=n+1,n$.
This is commutative for all $n$.

Since $\Lambda$ is a hull for $D_2$,
the compatible $\Lambda$-algebra structures on $A_{n+1},A_n,A_{n-1}$ induce
elements
$X_i'\in D_2(A_i)$, for $i=n-1,n$ and $n+1$. We may choose an object of
$\shD_2(A_{n+1})$ representing $X'_{n+1}$, which we also write as
$X'_{n+1}$, so that we have representatives $X_n'=\alpha_{n,*}X_{n+1}'$
and $X_{n-1}'=\alpha_{n-1,*}X_n'$. 

Let
$X_n\in D_0(A_n)\subseteq D_1(A_n)$ be arbitrary, so that $F(X_n)=X_n'$ as
elements of $D_2(A_n)$. Choosing a representative of the isomorphism class
of $X_n$, also written as $X_n$, we obtain an isomorphism 
$X_n'\mapright{\cong} F(X_n)$
in $\shD_2(A_n)$. Let $X_{n-1}=\alpha_{n-1,*}X_n$, 
$\psi_{n-1}:\epsilon_{n-1,*}X_n\rightarrow X_{n-1}$ the natural map
as in Theorem~\ref{1.8}. 
Let $v_{n-1}\in T^1_1(X_{n-1}/A_{n-1})$ be given by the pair
$(\epsilon_{n-1,*}X_n,\psi_{n-1})$. This gives 
$v'_{n-1}=F_*(v_{n-1})\in T^1_2(F(X_{n-1})/A_{n-1})$.

The isomorphism $X_{n-1}'\cong F(X_{n-1})$ induces a canonical isomorphism
\[
T^1_2(F(X_{n-1})/A_{n-1})\cong T^1_2(X_{n-1}'/A_{n-1})
\]
given by 
\[
(Y_{n-1}, Y_{n-1}\rightarrow F(X_{n-1}))\mapsto
(Y_{n-1}\times_{F(X_{n-1})} X_{n-1}',Y_{n-1}\times_{F(X_{n-1})} X_{n-1}'
\rightarrow X_{n-1}'),
\]
the fibre product existing by the (RS) condition. On the other hand,
by the co-fibred property,
we have a commutative diagram
\[
\xymatrix@C=20pt
{
X'_n\ar[r]^{\cong}\ar[d]& F(X_n)\ar[d]\\
\alpha_{n-1,*}X'_n\ar[r]_{\cong}& \alpha_{n-1,*}F(X_n)
}
\]
and hence a commutative diagram
\[
\xymatrix@C=20pt
{
\epsilon_{n-1,*}X'_n\ar[r]^{\cong}\ar[d]& \epsilon_{n-1,*}F(X_n)\ar[d]\\
\alpha_{n-1,*}X'_n\ar[r]_{\cong}& \alpha_{n-1,*}F(X_n)
}
\]
Using Definition~\ref{def:predef functor}, this diagram is the same
as
\[
\xymatrix@C=30pt
{
\epsilon_{n-1,*}X'_n\ar[r]\ar[d]&
F(\epsilon_{n-1,*}X_n)\ar[d]\\
X_{n-1}'\ar[r]& F(X_{n-1})
}
\]
Thus by Remark~\ref{rem:06J3}, this is a fibre diagram, and so
$v_{n-1}'\in T^1_2(F(X_{n-1}/A_{n-1}))$ is identified with
$v''_{n-1}=(\epsilon_{n-1,*}X'_n, \epsilon_{n-1,*}X'_n\rightarrow X'_{n-1})
\in T^1_2(X'_{n-1}/A_{n-1})$. Since $X_{n+1}'$ is a lifting of $X_n'$,
we see by Theorem \ref{1.8} that $v_{n-1}''$ is in the image of 
$T^1_2(X'_n/A_n)\rightarrow T^1_2(X'_{n-1}/A_{n-1})$, and hence
$v_{n-1}'$ is in the image of 
$T^1_2(F(X_n)/A_n)\rightarrow T^1_2(F(X_{n-1})/A_{n-1})$. 
Thus $\delta_2(v_{n-1}')=0$.  Since
$l|_{\im(\delta_1)}$ is injective, $\delta_1(v_{n-1})=0$. Again by Theorem \ref{1.8},
$X_n$ lifts to
$X_{n+1}\in D_1(A_{n+1})$, so $(d\otimes 1)\circ \delta_0(X_n)=\delta_1'(X_n)=0$, so
$\delta_0(X_n)\in\ker d\otimes 1= V\otimes (t^{n+1})$. Thus the image of $\delta_0$
is contained in
$V\otimes (t^{n+1})$. This proves the claim.
\end{proof}

We now follow some of the ideas of Kawamata in \cite{[19]} to construct $J'$. Dualizing
the sequence \eqref{(1.11)}, we have
$$J/\m_RJ\mapright{\phi} \m_{\Lambda}/\m^2_{\Lambda} \mapright{\dual{F_*}}
\m_S/\m_S^2.$$
Let $f_1,\ldots, f_t\in J$ be elements such that $\phi(f_1),\ldots,\phi(f_t)$
generate $\ker \dual{F_*}$, with $t=\dim_k \ker \dual{F_*}=\dim_k V$. Set $J'
=(f_1,\ldots,f_t)$. There is a surjective map
$\m_R/(\m_R^2+J')\rightarrow \m_R/(\m_R^2+J)$ and the dimensions of these spaces
are
the same by construction, so we have equality. 

The map 
$J'/\m_RJ'\mapright{i} J/\m_RJ$ induced by $J'\subseteq J$ is an inclusion since
\[
\dim_k J'/\m_RJ'\le t
\]
but $\dim_k \im(\phi\circ i)=t$. Thus
dually the composition
$$V\rightarrow \dual{(J/\m_RJ)} \rightarrow \dual{(J'/\m_RJ')}$$
is an isomorphism. Let $D_0':{\bf C}_{\Lambda}\rightarrow {\bf Ens}$ be the functor
pro-represented by $R/J'$. Clearly $D_0$ is a subfunctor of $D_0'$, and for each
$n$, we have a commutative diagram
\[
\xymatrix@C=30pt
{D_0(A_{n+1})\ar[r]\ar[d]&D_0(A_n)\ar[r]^{\delta_0}\ar[d]& V\otimes
(t^{n+1})\ar[d]^{\cong}\\
D_0'(A_{n+1})\ar[r]&D_0'(A_n)\ar[r]_{\delta_0'}&\dual{(J'/\m_RJ')}\otimes
(t^{n+1})
}
\]
This shows that an element of $D_0(A_n)$ lifts to $D_0(A_{n+1})$ if and only if
it lifts to $D_0'(A_{n+1})$. Since the choice of liftings is a principal
homogeneous space over $\dual{(\m_R/(\m_R\Lambda+\m_R^2))}\otimes (t^{n+1})$ for
both $D_0$ and $D_0'$, we see inductively that
$D_0(A_n) =D_0'(A_n)$ for any $\Lambda$-algebra structure on $A_n$. Thus
$$\hom_{\Lambda-alg}(R/J,k[[t]])=\hom_{\Lambda-alg}(R/J',k[[t]])$$
for any $\Lambda$-algebra structure on $k[[t]]$, and so
$$\hom_{k-alg}(R/J,k[[t]])=\hom_{k-alg}(R/J',k[[t]])$$
and we conclude that $\Supp(R/J)=\Supp(R/J')$ as in the last paragraph of the proof
of \cite[Thm.~1]{[19]}.\footnote{This is where the added hypothesis that
$k$ is algebraically closed is necessary.} 
\end{proof}

\section{Obstructions for Calabi-Yau threefolds with canonical singularities}
\label{sec2}

\begin{say}
\label{(2.1)}
Let $X$ be a Calabi-Yau threefold with canonical singularities over $k=\CC$,
the complex numbers, i.e., a projective threefold with canonical singularities,
$K_{X}=0$, and $h^1(\O_{X})=0$. Set $Z=\Sing(X)$. Note $Z$ can be zero or one dimensional.
We want to study the deformation
theory of
$X$. In particular, we will relate the deformation theory of $X$ to the deformation
theory of $\hat X$, the formal completion of $X$ along $Z$. To paraphrase Theorem
\ref{2.2} below, we will find that the obstructions to deforming $X$ are contained in the
obstructions to deforming $\hat X$, i.e. ``obstructions to deforming $X$ are local
to the singularities of $X$.''

To make this concept rigorous, let $D$ be the functor of deformations of $X$, 
and let $D_{\loc}$ be the functor of deformations of $\hat X$, as in 
Example~\ref{ex:principal}. There
is a natural morphism of functors $F:D\rightarrow D_{\loc}$ taking a deformation
to its completion along $Z$. Note that $D$ is pro-representable
since
$\hom(\Omega^1_X,\O_X)=0$ by \cite[Cor.~8.6]{[17]}.

If $X_n$ is a deformation of
$X$ over $A_n$, we have
$T^1(X_n/A_n)\cong
\ext^1_{\O_{X_n}}(\Omega^1_{X_n/A_n},\O_{X_n})$ and $T^1_{\loc}(X_n/A_n)\cong
\ext^1_{\O_{\hat X_n}}(\Omega^1_{X_n/A_n}\otimes_{\O_{X_n}} \O_{\hat X_n},
\O_{\hat
X_n})$. The former isomorphism is well-known. For the latter, if $X_n/A_n$
is locally embedded in $Y_n/A_n$ smooth with ideal sheaf $\I$, then the local
${\bf T}^1$ sheaf of $\hat X_n/A_n$ is as usual given by
$$\lhom_{\O_{\hat X_n}}(\Omega^1_{\hat Y_n/A_n}|_{\hat X_n}, \O_{\hat X_n})
\rightarrow \lhom_{\O_{\hat X_n}}(\widehat{\I/\I^2}, \O_{\hat X_n})\rightarrow
{\bf T}^1\rightarrow 0.$$
Note that $\Omega^1_{Y_n/A_n}\otimes_{\O_{Y_n}} \O_{\hat Y_n} \cong
\hat\Omega^1_{\hat Y_n/A_n}$, the completion of $\Omega^1_{\hat Y_n/A_n}$
(see \cite[Chap.~0, 20.7.14]{[14]}). Since $\O_{\hat X_n}$ is complete, 
\begin{align*}
\lhom_{\O_{\hat X_n}}(\Omega^1_{\hat Y_n/A_n}|_{\hat X_n}, \O_{\hat
X_n})= {} &
\lhom_{\O_{\hat X_n}}(\hat\Omega^1_{\hat Y_n/A_n}|_{\hat X_n}, \O_{\hat X_n})\\
={} &
\lhom_{\O_{\hat X_n}}(\Omega^1_{Y_n/A_n}|_{X_n}\otimes_{\O_{X_n}}
\O_{\hat X_n}, \O_{\hat X_n}).
\end{align*}
Thus from the exact sequence
$$\I/\I^2\otimes_{\O_{X_n}} \O_{\hat X_n}\rightarrow \Omega^1_{Y_n/A_n}|_{X_n}\otimes_{\O_{X_n}}
\O_{\hat X_n} \rightarrow \Omega^1_{X_n/A_n}\otimes_{\O_{X_n}} \O_{\hat X_n}
\rightarrow 0$$
we see that ${\bf T}^1\cong \lext^1_{\O_{\hat X_n}}(\Omega^1_{X_n/A_n}
\otimes_{\O_{X_n}} \O_{\hat X_n},
\O_{\hat
X_n})$. Local infinitesimal deformations of $\hat X_n$ then patch together as usual
to yield an element of
$\ext^1_{\O_{\hat X_n}}(\Omega^1_{X_n/A_n}\otimes_{\O_{X_n}} \O_{\hat X_n},
\O_{\hat
X_n})$.

Finally, let $T^2=\ext^2_{\O_X}(\Omega^1_X, \O_X)$ and $T^2_{\loc}=\ext^2_{\O_{X}}
(\Omega^1_X,\O_{\hat X})$. Let $l:T^2\rightarrow T^2_{\loc}$ be the map induced
by the map $\O_X\rightarrow\O_{\hat X}$. 
\end{say}

The following theorem is our generalization of the Bogomolov-Tian-Todorov
unobstructedness theorem.

\begin{theorem}
\label{2.2} 
Let $X$ be a Calabi-Yau threefold with canonical
singularities, and $X_n/A_n$ a deformation of $X$. There is a commutative diagram
\[
\xymatrix@C=30pt
{
T^1(X_n/A_n)\ar[r]^{F}\ar[d]&T^1_{\loc}(X_n/A_n)\ar[d]\\
T^1(X_{n-1}/A_{n-1})\ar[r]^{F}\ar[d]_{\delta}
&T^1_{\loc}(X_{n-1}/A_{n-1}) \ar[d]^{\delta_{\loc}}
\\
T^2\ar[r]_{l}&T^2_{\loc}
}
\]
where
$l|_{\im(\delta)}:\im(\delta)\rightarrow T^2_{\loc}$ is injective.
\end{theorem}

We shall prove this result by generalising Namikawa's argument in \cite{[28]}.

\begin{lemma}
\label{2.3} 
Given the hypotheses of \ref{(2.1)}, if $X_n$ is a deformation
of $X$ over $A_n$, then there are natural isomorphisms 
$$\ext^i_{\O_{X_n}}(\Omega^1_{X_n/A_n},\O_{X_m})\cong
\ext^i_{\O_{X_m}}(\Omega^1_{X_m/A_m},\O_{X_m})$$
for $m<n, i\le 2$ and $X_m=X_n\otimes_{A_n} A_m$.
\end{lemma}

\begin{proof}
The change of rings spectral sequence (\cite[Thm.~10.74]{[40]}) tells us that
$$\lext^p_{\O_{X_m}}({\rm Tor}^{\O_{X_n}}_q(\Omega^1_{X_n/A_n},\O_{X_m}),\O_{X_m})
\Rightarrow \lext^i_{\O_{X_n}}(\Omega^1_{X_n/A_n},\O_{X_m}).$$
Now $\Omega^1_{X_n/A_n}$ is a flat $\O_{X_n}$-module away from $Z\subseteq X_n$, so
${\rm Tor}^{\O_{X_n}}_q(\Omega^1_{X_n/A_n},\O_{X_m})$ is supported on $Z$ for $q\ge 1$.
Since $X_m$ is Cohen-Macaulay,
\[
\lext^p_{\O_{X_m}}({\rm Tor}^{\O_{X_n}}_q(\Omega^1_{X_n/A_n},\O_{X_m}),\O_{X_m})=0
\]
for $p\le 1, q\ge 1$. Thus
\begin{align*}
\lext^i_{\O_{X_m}}(\Omega^1_{X_m/A_m},\O_{X_m})= {} &
\lext^i_{\O_{X_m}}(\Omega^1_{X_n/A_n}\otimes_{\O_{X_n}}\O_{X_m},\O_{X_m})\\
={} & \lext^i_{\O_{X_n}}(\Omega^1_{X_n/A_n},\O_{X_m})
\end{align*}
for $i\le 2$. The statement of the lemma for global Exts then follows from the
local-global spectral sequence for Exts.
\end{proof}

\begin{lemma}
\label{2.4} 
There are isomorphisms
\begin{itemize}
\item[(a)] $\ext^1_{\O_{X_n}}(\Omega^1_{X_n/A_n},\O_{\hat X_m})
\cong T^1_{\loc}(X_m/A_m)$ for $m<n$. 
\item[(b)] $\ext^2_{\O_{X_n}}(\Omega^1_{X_n/A_n}, \O_{\hat X}) = T^2_{\loc}.$
\item[(c)] $\dual{(T^2_{\loc})}\cong H^1_Z(\Omega^1_X)$, $\dual{(T^2)}
\cong H^1(\Omega^1_X)$, and the natural map $H^1_Z(\Omega^1_X)\rightarrow
H^1(\Omega^1_X)$ is dual to $l$.
\end{itemize}
\end{lemma}

\begin{proof}
The change of rings spectral sequence yields
$$\lext^p_{\O_{X_m}}({\rm Tor}^{\O_{X_n}}_q(\Omega^1_{X_n/A_n},\O_{X_m}),\O_{\hat X_m})
\Rightarrow \lext^i_{\O_{X_n}}(\Omega^1_{X_n/A_n},\O_{\hat X_m}),$$
so the same argument as in the proof of Lemma \ref{2.3} shows that
$$\ext^i_{\O_{X_n}}(\Omega^1_{X_n/A_n},\O_{\hat X_m})\cong
\ext^i_{\O_{X_m}}(\Omega^1_{X_m/A_m},\O_{\hat X_m})$$
for $i\le 2$. Furthermore, since $\O_{\hat X_m}$
is a flat $\O_{X_m}$-module,
$$\ext^i_{\O_{X_m}}(\Omega^1_{X_m/A_m},\O_{\hat X_m})
\cong\ext^i_{\O_{\hat X_m}}(\Omega^1_{X_m/A_m}\otimes_{\O_{X_m}} \O_{\hat X_m},
\O_{\hat X_m}),$$ which is $T^1_{\loc}(X_m/A_m)$ by \ref{(2.1)} for $i=1$.
For $i=2$, $m=0$, this gives (b).

(c) follows from the usual statement of local
duality if $Z$ is dimension zero. Since $X$ may not have isolated singularities, 
we use Alonso, Jerem\'\i as and Lipman's generalization of local
duality \cite{[1]}. By \cite[(0.3), (0.1), and (0.4.1)]{[1]}, 
$${\bf R}\lhom^{\bullet}_{\O_X}({\bf R}{\underline\Gamma}_Z\Omega^1_X,\O_X)
\cong{\bf R}\lhom^{\bullet}_{\O_X}(\Omega^1_X,\O_{\hat X}).$$
By the Grothendieck duality theorem of \cite{[15]}, 
\begin{align*}
{\bf R}\Gamma{\bf R}\lhom^{\bullet}_{\O_X}({\bf
R}{\underline\Gamma}_Z\Omega^1_X,\O_X)\cong {} & \hom^{\bullet}_k({\bf R}\Gamma{\bf
R}{\underline\Gamma}_Z\Omega^1_X,k[-3])\\
= {} &\hom^{\bullet}_k({\bf R}\Gamma_Z\Omega^1_X,k[-3]).
\end{align*}
Putting these two isomorphisms together and taking cohomology of the two complexes,
we have
$$H^i(\hom^{\bullet}_k({\bf R}\Gamma_Z\Omega^1_X,k[-3]))\cong H^i({\bf
R}\hom^{\bullet}_{\O_X}(\Omega^1_X,\O_{\hat X})).$$
For $i=2$, this yields the isomorphism
$$\dual{(H^1_Z(\Omega^1_X))}\cong \ext^2_{\O_X}(\Omega^1_X,\O_{\hat X})\cong
T^2_{\loc}.$$
This isomorphism is compatible with the map
\[
T^2=\ext^2_{\O_X}(\Omega^1_X,\O_X)\rightarrow
\ext^2_{\O_X}(\Omega^1_X,\O_{\hat X})=T^2_{\loc}
\]
induced by the map $\O_X\rightarrow
\O_{\hat X}$ since \cite[(0.3)]{[1]} also tells us that the following diagram is
commutative:
\[
\xymatrix@C=30pt
{
{\bf R}\lhom_{\O_X}^{\bullet}(\Omega^1_X,\O_X)
\ar[d]_{\alpha}\ar[rd]^{\beta}
&\\
{\bf R}\lhom^{\bullet}_{\O_X}({\bf R}{\underline\Gamma}_Z\Omega^1_X,\O_X)
\ar[r]_{\cong}&{\bf R}\lhom^{\bullet}_{\O_X}(\Omega^1_X,\O_{\hat X})
}
\]
where $\alpha$ is induced by the natural map in the derived category ${\bf
R}{\underline\Gamma}_Z\Omega^1_X
\rightarrow \Omega^1_X$ and $\beta$ is induced by $\O_X\rightarrow \O_{\hat X}$.
\end{proof}

We need a version of Namikawa's \cite[Lem.~2.2]{[28]}.
Let $U=X-Z$. Then
$X_n/A_n$ induces a deformation $U_n$ of $U$ over $A_n$. We set $\tilde
\Omega^1_{X_n/A_n}= j_*\Omega^1_{U_n/A_n}$ with $j:U_n\rightarrow X_n$ the
inclusion. We put $\tilde \Omega^1_X:=\tilde \Omega^1_{X/k}$. 
Observe that we have a natural map
$d\log:H^1(X,\O_X^*)\rightarrow H^1(X,\tilde \Omega^1_X)$ which is
the composition of $H^1(X,\O_X^*)\rightarrow H^1(X,\Omega^1_X)$ and 
$H^1(X,\Omega^1_X)\rightarrow H^1(X,\tilde \Omega^1_X)$. Also, if $\pi:\tilde X
\rightarrow X$ is a resolution of singularities of $X$, then the
map
$\Omega_X^1\rightarrow \tilde\Omega_X^1$ factors through $\Omega_X^1
\rightarrow \pi_*\Omega^1_{\tilde X}$.

\begin{lemma}
\label{2.5} 
The image of the map 
$$d\log:H^1(X,\O_X^*)\rightarrow H^1(X,\tilde\Omega^1_X)$$
generates $H^1(X,\tilde\Omega_X^1)$ as a $k$-vector space. 
\end{lemma}

\begin{proof}
The proof of \cite[Lem.~2.2]{[28]} actually shows in general that
the image of $H^1(X,\O_X^*)\rightarrow H^1(X,\pi_*\Omega^1_{\tilde X})$
generates the latter as a vector space, using only the hypothesis
that $\tilde X$ has rational singularities, which is true of any
canonical singularity. Now the kernel and cokernel of
$\pi_*\Omega^1_{\tilde X}\rightarrow \tilde\Omega^1_X$ are supported on $Z$.
If we can show furthermore that $\coker(\pi_*\Omega^1_{\tilde X}\rightarrow
\tilde \Omega^1_X)$ is supported on a finite set of points, then in fact
$H^1(\pi_*\Omega^1_{\tilde X})\rightarrow H^1(\tilde\Omega^1_X)$ is surjective
and the lemma follows. 

By \cite{[38]}, except for a finite number of dissident points,
$X$ is analytically isomorphic in a neighborhood of a point of $Z$ to 
$\Delta\times S$, where $\Delta$ is the germ of a curve and $S$ is a germ of
a du Val surface singularity, with resolution $\pi':\tilde S\rightarrow S$. Thus, in 
this neighborhood, $\pi:\tilde X\rightarrow X$ looks like $\pi:\Delta\times\tilde S
\rightarrow \Delta\times S$. Let $p_1$ and $p_2$ be the projections of
$\Delta\times S$ onto $\Delta$ and $S$ respectively. The map
$\pi_*\Omega^1_{\tilde X}\rightarrow
\tilde \Omega^1_X$ is $p_1^*\Omega^1_{\Delta}\oplus p_2^*\pi'_*\Omega^1_{\tilde
S}\rightarrow
p_1^*\Omega^1_{\Delta}\oplus p_2^*\tilde\Omega^1_S$. Since $S$ is a rational
singularity,  it follows from \cite{[43]} that $\pi_*\Omega^1_{\tilde S}\rightarrow \tilde
\Omega^1_S$ is surjective. Thus $\pi_*\Omega^1_{\tilde
X}\rightarrow\tilde\Omega^1_X$ is surjective in this neighborhood, and so
$\coker(\pi_*\Omega^1_{\tilde X}\rightarrow\tilde\Omega^1_X)$ is concentrated on the
dissident points of $Z$.
\end{proof}

\begin{proof}[Proof of the theorem]
 To obtain the first column of the diagram, we apply the functor
$\hom_{\O_{X_n}}(\Omega^1_{X_n/A_n},\cdot)$ to the exact sequence
$$\exact{\O_X}{\O_{X_n}}{\O_{X_{n-1}}}$$
yielding
\begin{equation}
\label{eq:ext sequence}
\ext^1_{\O_{X_n}}(\Omega^1_{X_n/A_n},\O_{X_n})\rightarrow
\ext^1_{\O_{X_n}}(\Omega^1_{X_n/A_n},\O_{X_{n-1}})\rightarrow
\ext^2_{\O_{X_n}}(\Omega^1_{X_n/A_n},\O_X),
\end{equation}
from which we obtain the first column, by Lemma \ref{2.3}.
The second column is obtained by
applying
$\hom_{\O_{X_n}}(\Omega^1_{X_n/A_n},\cdot)$ to the sequence
$$\exact{\O_{\hat X}}{\O_{\hat X_n}}{\O_{\hat X_{n-1}}},$$ and applying
Lemma \ref{2.4}. The maps between the columns induced by the maps $\O_{X_n}\rightarrow
\O_{\hat X_n}$
yield a commutative diagram. For the first two rows, these maps coincide with
the maps induced by $F:D\rightarrow D_{\loc}$.

We now only need to
show that $l|_{\im(\delta)}:\im(\delta)\rightarrow T^2_{\loc}$ is injective for each
such diagram. This is equivalent to $\dual{l|_{\im(\delta)}}:\dual{(T^2_{\loc})}
\rightarrow \dual{\im(\delta)}$ being surjective. 

We have a diagram
\begin{equation}
\label{eq:big diagram}
\xymatrix@=30pt
{H^1_Z(\Omega^1_{X_n/A_n})\ar[r]\ar[d]&H^1(\Omega^1_{X_n/A_n})
\ar[r]\ar[d]&H^1(\tilde\Omega^1_{X_n/A_n})\ar[d]&\\
H^1_Z(\Omega^1_X)\ar[r]\ar[d]&H^1(\Omega^1_X)\ar[d]
\ar[r]&\ar[d]H^1(\tilde\Omega^1_X)\ar[d]&\\
K_{\loc}\ar[r]\ar[d]&K\ar[r]\ar[d]&0\\
0&0&
}
\end{equation}
Here the vertical
maps between the first two rows are induced by the restriction map
$\Omega^1_{X_n/A_n}\rightarrow \Omega^1_X$, and $K_{\loc}$ and $K$ are
just defined to make the columns exact. The map $H^1_Z(\Omega^1_X)
\rightarrow H^1(\Omega^1_X)$ is dual to $l:T^2\rightarrow T^2_{\loc}$, by
Lemma \ref{2.4} (c). I claim the rows of this diagram are exact. 

{\it Exactness of the first two rows:} For any $n$, we have the exact sequence
$$0\rightarrow \HH^0_Z(\Omega^1_{X_n/A_n})\rightarrow \Omega^1_{X_n/A_n}
\rightarrow \tilde\Omega^1_{X_n/A_n} \rightarrow \HH^1_Z(\Omega^1_{X_n/A_n})
\rightarrow 0.$$
Now $\HH^0_Z(\Omega^1_{X_n/A_n})$ is supported on the set of points of $X$ which 
are not locally complete intersection, and so has finite support. Thus
$H^1(\HH^0_Z(\Omega^1_{X_n/A_n}))=0$ and $H^1_Z(\Omega^1_{X_n/A_n})
=H^0(\HH^1_Z(\Omega^1_{X_n/A_n}))$. From this
follows exactness of the first two rows. 

{\it Exactness of the last row:} As in
\cite{[28]}, there is a diagram
\[
\xymatrix@C=30pt
{
H^1(\O_{X_n}^*)\ar[d]_{\alpha}\ar[r]&H^1(\tilde\Omega^1_{X_n/A_n})\ar[d]\\
H^1(\O_X^*)\ar[r]&H^1(\tilde\Omega^1_X)
}
\]
with $\alpha$ surjective since $H^2(\O_X)=0$, and by Lemma \ref{2.5}, the image
of $H^1(\O_X^*)$ in $H^1(\tilde \Omega^1_X)$ generates the latter as a
$k$-vector space. Thus the composed map $H^1(\O_{X_n}^*)\otimes_{\ZZ} k
\rightarrow H^1(\Omega^1_{X_n/A_n})\rightarrow H^1(\tilde\Omega^1_{X_n/A_n})
\rightarrow H^1(\tilde\Omega^1_X)$ is surjective, and so the map 
$H^1(\Omega^1_{X_n/A_n})\rightarrow H^1(\tilde\Omega^1_X)$ is surjective.
Now a simple diagram chase shows that the last row is exact.

\medskip

The diagram \eqref{eq:big diagram} now shows that 
the composed map $H^1_Z(\Omega^1_X)\rightarrow K$ is surjective. 

We now consider the following isomorphisms of $\partial$-functors (i.e., these
functors take triangles to triangles and the isomorphisms preserve triangles):
\begin{align}
\label{key iso}
\begin{split}
{\bf R}\Hom_{A_n}(Rf_*(\Omega_{X_n/A_n}\otimes^L \bullet),A_n)
\cong {} &  {\bf R}\Hom_{\O_{X_n}}(\Omega_{X_n/A_n}\otimes^L \bullet, \omega_{X_n/A_n}[3])\\
\cong {} & {\bf R}\Hom_{\O_{X_n}}(\Omega_{X_n/A_n},
{\bf R}\lhom_{\O_{X_n}}(\bullet,\omega_{X_n/A_n}[3])),
\end{split}
\end{align}
where $f:X_n\rightarrow \Spec A_n$ is the structure morphism.
Here the first isomorphism is Grothendieck duality as stated in
\cite[Thm.~(4.1.1)]{Lip}, which also states it is a $\partial$-functor.
The second isomorphism is \cite[Prop.~5.15]{[15]} or \cite[Prop.~(2.6.1)]{Lip}.
Following \cite{[15]}, the first isomorphism
requires $\Omega_{X_n/A_n}\otimes^L\bullet$ to lie in $D^-_{qc}(X_n)$,
i.e., complexes of quasicoherent sheaves bounded above,
and the second isomorphism requires $\bullet$ to lie in
$D^-_c(X)$, i.e., complexes of coherent sheaves bounded above.
We wish to apply this for $\bullet=\O_{X_m}$ for various $m\le n$.
Since this sheaf has a resolution by locally free sheaves
\begin{equation}
\label{eq:resolution}
\xymatrix@C=30pt
{\cdots\ar[r]&\O_{X_n}\ar[r]^{t^{m+1}}&\O_{X_n}\ar[r]^{t^{n-m}}&
\O_{X_n}\ar[r]^{t^{m+1}}&\O_{X_n}\ar[r]&\O_{X_{m}}\ar[r]&0,}
\end{equation}
$\Omega_{X_n/A_n}\otimes^L \O_{X_m}$ is represented by a complex
bounded above, as needed. 

We apply these isomorphisms to the triangle arising from the exact sequence
\begin{equation}
\label{eq:triangle}
0\rightarrow \O_{X_{n-1}} \rightarrow \O_{X_n}\rightarrow \O_X\rightarrow 0.
\end{equation}
Because these isomorphisms are isomorphisms of $\partial$-functors,
after taking cohomology of the left and right sides of \eqref{key iso},
we obtain isomorphic long exact sequences of
$A_n$-modules, the key part of these being, for us:
\begin{align}
\label{eq:exact1}
\begin{split}
&H^{-2}( {\bf R}\Hom_{A_n}(Rf_*(\Omega_{X_n/A_n}\otimes^L \O_{X_n}),A_n))\\
\rightarrow {} & 
H^{-2}({\bf R}\Hom_{A_n}(Rf_*(\Omega_{X_n/A_n}\otimes^L \O_{X_{n-1}}),A_n))\\
\rightarrow {} &H^{-1}({\bf R}\Hom_{A_n}(Rf_*(\Omega_{X_n/A_n}\otimes^L \O_X),A_n))
\end{split}
\end{align}
and 
\begin{align}
\label{eq:exact2}
\begin{split}
&H^{-2}({\bf R}\Hom_{\O_{X_n}}(\Omega_{X_n/A_n},
{\bf R}\lhom_{\O_{X_n}}(\O_{X_n},\omega_{X_n/A_n}[3])))\\
\rightarrow {} & 
H^{-2}({\bf R}\Hom_{\O_{X_n}}(\Omega_{X_n/A_n},
{\bf R}\lhom_{\O_{X_n}}(\O_{X_{n-1}},\omega_{X_n/A_n}[3])))\\
\rightarrow {} & H^{-1}({\bf R}\Hom_{\O_{X_n}}(\Omega_{X_n/A_n},
{\bf R}\lhom_{\O_{X_n}}(\O_X,\omega_{X_n/A_n}[3])))
\end{split}
\end{align}

Let us first compute the terms in \eqref{eq:exact1}. As $A_n$ is
an injective $A_n$-module, we may replace ${\bf R}\Hom_{A_n}(\bullet,
A_n)$ with $\Hom_{A_n}(\bullet,A_n)$. Thus \eqref{eq:exact1}
may be written as
\begin{align*}
\Hom_{A_n}(\HH^2(X_n,\Omega_{X_n/A_n}\otimes^L \O_{X_{n}}),A_n)
& \rightarrow
\Hom_{A_n}(\HH^2(X_n,\Omega_{X_n/A_n}\otimes^L \O_{X_{n-1}}),A_n)\\ &
\rightarrow \Hom_{A_n}(\HH^1(X_n,\Omega_{X_n/A_n}\otimes^L \O_X),A_n).
\end{align*}
As $\Omega_{X_n/A_n}\otimes^L \O_{X_n}\cong \Omega_{X_n/A_n}$, the first
term is just $\Hom_{A_n}(H^2(X_n,\Omega_{X_n/A_n}),A_n)$. The other
two terms can be computed via the hypercohomology spectral sequence
(for $m=n-1$ or $0$) given at the $E^2$ term by
\[
H^p(X_n,\shH^q(\Omega_{X_n/A_n}\otimes^L\O_{X_m}))
\Rightarrow \HH^{p+q}(X_n,\Omega_{X_n/A_n}\otimes^L\O_{X_m}),
\]
where $\shH^*$ denotes cohomology of a complex of sheaves.
Now 
\[
\shH^q(\Omega_{X_n/A_n}\otimes^L\O_{X_m})=\Tor_{-q}^{\O_{X_n}}
(\Omega_{X_n/A_n},\O_{X_m}).
\]
The latter is $\Omega_{X_m/A_m}$ for
$q=0$ and is supported on $\Sing(X)$ for $q<0$. Thus
the spectral sequence is fourth quadrant and
has the following terms, the upper left entry
being the $(0,0)$ entry:
\[
\begin{matrix}
H^0(X_n,\Omega_{X_m/A_m})&
H^1(X_n,\Omega_{X_m/A_m})&
H^2(X_n,\Omega_{X_m/A_m})&
H^3(X_n,\Omega_{X_m/A_m})\\
?&?&0&0\\
?&?&0&0
\end{matrix}
\]
The spectral sequence thus degenerates at the $E_2$ term, and one
sees that
\[
\HH^i(X_n,\Omega_{X_n/A_n}\otimes^L\O_{X_m})
\cong H^i(X_n,\Omega_{X_m/A_m})=H^i(X_m,\Omega_{X_m/A_m})
\]
for $i=1,2$.
In particular, the long exact sequence of hypercohomology obtained by
applying the functor $\HH^{\bullet}(X_n,\Omega_{X_n/A_n}\otimes^L \bullet)$
to the triangle obtained from \eqref{eq:triangle} produces the long
exact sequence
\begin{equation}
\label{eq:second long}
H^1(\Omega_{X_n/A_n})\mapright{} H^1(\Omega_X)
\mapright{} H^2(\Omega_{X_{n-1}/A_{n-1}}) \mapright{\cdot t}
H^2(\Omega_{X_n/A_n})
\end{equation}
showing that 
\[
K=\ker(H^2(\Omega_{X_{n-1}/A_{n-1}})
\mapright{\cdot t} H^2(\Omega_{X_n/A_n}))
=\im(H^1(\Omega_X)\rightarrow H^2(\Omega_{X_n/A_n})).
\]
Also, noting that $\Hom_{A_n}(M,A_n)=\Hom_{A_m}(M,A_m)$ for $M$ an $A_m$-module,
we get that \eqref{eq:exact1} coincides with
\begin{equation}
\label{eq:dual sequence}
\Hom_{A_n}(H^2(\Omega_{X_n/A_n}),A_n)
\rightarrow
\Hom_{A_{n-1}}(H^2(\Omega_{X_{n-1}/A_{n-1}}),A_{n-1})
\rightarrow
\Hom_k(H^1(\Omega_{X/k}),k),
\end{equation}
the dual of the last three terms of \eqref{eq:second long}.

We turn to \eqref{eq:exact2}. We 
calculate for $m\le n$, and $i=-2,-1$:
\begin{align*}
&H^{i}({\bf R}\Hom_{\O_{X_n}}(\Omega_{X_n/A_n},
{\bf R}\lhom_{\O_{X_n}}(\O_{X_m},\omega_{X_n/A_n}[3])))\\
\cong {} &
H^{i+3}({\bf R}\Hom_{\O_{X_n}}(\Omega_{X_n/A_n},
{\bf R}\lhom_{\O_{X_n}}(\O_{X_m},\omega_{X_n/A_n})))\\
\cong {}  & \hyperext^{i+3}_{\O_{X_n}}(\Omega_{X_n/A_n},{\bf R}\lhom
(\O_{X_m},\omega_{X_n/A_n}))\\
\cong {} & \hyperext^{i+3}_{\O_{X_n}}(\Omega_{X_n/A_n},{\bf R}\lhom
(\O_{X_m},\O_{X_n}))
\end{align*}
using the Calabi-Yau condition and the hyperext.

However, direct calculation from the resolution \eqref{eq:resolution}
of $\O_{X_m}$ shows that
\[
\Ext^q_{\O_{X_n}}(\O_{X_m},\O_{X_n})\cong 
\begin{cases}
\O_{X_m} & q=0;\\
0 & q>0.
\end{cases}
\]
Thus in fact
\begin{align*}
&H^{i}({\bf R}\Hom_{\O_{X_n}}(\Omega_{X_n/A_n},
{\bf R}\lhom_{\O_{X_n}}(\O_{X_m},\omega_{X_n/A_n}[3])))\\
\cong {} & \Ext^{i+3}_{\O_{X_n}}(\Omega_{X_n/A_n},\O_{X_m})\\
\cong {} & \Ext^{i+3}_{\O_{X_m}}(\Omega_{X_m/A_m},\O_{X_m})
\end{align*}
for $i+3\le 2$ by Lemma \ref{2.3}. This shows that the right-hand
vertical side of the diagram of interest agrees with
\begin{equation}
\label{eq:exts}
\Ext^1_{\O_{X_n}}(\Omega_{X_n/A_n},\O_{X_n})
\mapright{}
\Ext^1_{\O_{X_{n-1}}}(\Omega_{X_{n-1}/A_{n-1}},\O_{X_{n-1}})
\mapright{\delta}
\Ext^2_{\O_X}(\Omega_{X},\O_{X}).
\end{equation}
Putting this together with \eqref{eq:dual sequence}, we see we obtain a commutative
diagram with vertical arrows being isomorphisms and rows being exact:
\begin{equation}
\label{eq:bob diagram}
\xymatrix@C=30pt
{
\Hom_{A_n}(H^2(\Omega_{X_n/A_n}),A_n)\ar[r]\ar[d]&
\Hom_{A_{n-1}}(H^2(\Omega_{X_{n-1}/A_{n-1}}),A_{n-1})
\ar[r]\ar[d]&
\Hom_k(H^1(\Omega_{X/k}),k)\ar[d]\\
\Ext^1_{\O_{X_n}}(\Omega_{X_n/A_n},\O_{X_n})\ar[r]&
\Ext^1_{\O_{X_{n-1}}}(\Omega_{X_{n-1}/A_{n-1}},\O_{X_{n-1}})
\ar[r]_>>>>>>{\delta}&
\Ext^2_{\O_X}(\Omega_{X},\O_{X})
}
\end{equation}
This shows that $K$ coincides with $\dual{(\im\delta)}$.

So the map $H^1(\Omega^1_X)\rightarrow K\rightarrow 0$ is dual to $0\rightarrow
\im\delta\rightarrow T^2$. Thus
$\dual{(T^2_{\loc})}\rightarrow \dual{\im(\delta)}$ is the surjection
$H^1_Z(\Omega^1_X)\rightarrow K$. This is the desired surjectivity. 
\end{proof}

\begin{remark}
\label{2.6}
We cannot apply Theorem \ref{2.2} immediately to the situation of
Theorem \ref{1.9} without first knowing that $D_{\loc}$ has a hull. By 
\ref{(1.5)}, this is
the case if and only if $T^1_{\loc}(X/k)$ is finite dimensional. This is of course
the case if $X$ has isolated singularities, in which case, we can just as well
consider the complex germ $(X,Z)$ instead of the formal scheme $\hat X$.
However, $T^1_{\loc}$ need not be finite dimensional
if $X$ has non-isolated
singularities. Nevertheless, even in this case, Theorem \ref{2.2} can be useful. In any
event, Theorem \ref{2.2} tells us that 
$$T^1(X_n/A_n)\rightarrow T^1(X_{n-1}/A_{n-1})\rightarrow T^2_{\loc}$$
is exact. If $T^2_{\loc}=0$, then $D$ is unobstructed by the $T^1$-lifting
criterion. If $T^2_{\loc}\not=0$, then we obtain dimension estimates for
$\Def(X)$ using methods similar to \cite{[19]}. See \cite{[12]} 
for an application in the
non-isolated case. Note in particular that if $X$ has isolated complete
intersection singularities, then $T^2_{\loc}=0$, reproducing the unobstructedness
result of \cite{[28]} in a rather more inefficient way.
\end{remark}

\begin{remark}
While the proof of Theorem \ref{2.2} only works in dimension three,
the higher-dimensional analogue of the diagram \eqref{eq:bob diagram}
still exists for all $\dim X\ge 3$, with $H^2$ and $H^1$ replaced by $H^{\dim X-1}$ and
$H^{\dim X-2}$ respectively. The proof is the same.
\end{remark}

\begin{example}
\label{2.7}
We give here a simple example of a Calabi-Yau with obstructed
deformation theory. Let $Y\subseteq\P^9$ be a cone over a non-singular del Pezzo
surface in $\P^8$ isomorphic to $\Ptwo$ blown up in one point. By \cite{[2]}, 
a hull of the deformation functor of the singular point of $Y$ is
$\Lambda=k[t]/(t^2)$. In fact, in \cite[(9.2)]{[2]}, Altmann gives explicit equations
for a deformation of the affine cone in ${\bf C}^9$ over the del Pezzo surface
in $\P^8$, and these equations are easily projectivized to yield a non-trivial but
obstructed infinitesimal deformation of $Y$, $\shY/\Spec\Lambda$. 

Let $X$ be a double cover of $Y$ branched over the intersection of $Y$ with a
general quartic hypersurface in $\P^9$. It is easy to see that $K_X=0$, $X$ has two
singular points analytically isomorphic to the singular point of $Y$, and
the deformation $\shY/\Spec\Lambda$ lifts to a deformation $\X/\Spec\Lambda$, 
which is obstructed. In fact, $\Def(X)$ is non-reduced.

See \cite{[11]} for an example in any dimension of obstructed deformations for a
Calabi-Yau with non-isolated singularities. Such an example is much more subtle,
and requires a more global analysis. 
\end{example}

\section{Calabi-Yaus with complete intersection singularities}
\label{sec3}

\begin{say}
\label{(3.1)} We will first consider, very generally, the situation that $(X,0)$ is the germ
of an isolated rational complex threefold singularity, and that $\pi:(\tilde
X,E)\rightarrow (X,0)$ is a resolution of singularities. We have a natural map of
germs of analytic spaces
$\Def(\tilde X)\rightarrow \Def(X)$ by \cite[Pro.~11.4]{[23]}, 
(by \cite[Thm.~1.4 (c)]{[48]}
for the map on the level of deformation functors,)
since $H^1(\O_{\tilde X})=0$. We denote by $\O_{X,0}$
the local ring of $X$ at the origin with maximal ideal $\m$, and we denote by $T^1$
the tangent space of
$\Def(X)$.
\end{say}

\begin{lemma}
\label{3.2} 
The tangent space to $\Def(\tilde X)$ is $H^0(R^1\pi_*
\T_{\tilde X})$, of $\Def(X)$ is $T^1=H^2_Z(\T_X)$, where $\T_X=\hom_{\O_X}
(\Omega^1_X,\O_X)$, and there is an exact sequence of $\O_{X,0}$-modules
$$H^0(R^1\pi_*\T_{\tilde X})\rightarrow H^2_Z(\T_X) \rightarrow T'
\rightarrow 0$$ with $T'=\ker(
H^2_E(\T_{\tilde X})\rightarrow H^0(R^2\pi_*\T_{\tilde X}))$,
$Z=\Sing(X)$ and $E$ the exceptional locus of $\pi$. The map
$H^0(R^1\pi_*\T_{\tilde X})\rightarrow H^2_Z(\T_X)$ is the differential
of the map $\Def(\tilde X)\rightarrow \Def(X)$. 
\end{lemma}

\begin{proof}
Since $X$ is a germ, $H^1(\pi_*\T_{\tilde X})=H^2(\pi_*\T_{\tilde X})=0$.
Thus the tangent space to $\Def(\tilde X)$ is
$H^1(\T_{\tilde X})\cong H^0(R^1\pi_*\T_{\tilde X})$ by the Leray spectral
sequence. Similarly, $H^2(\T_{\tilde X})=H^0(R^2\pi_*\T_{\tilde X})$. Also,
the tangent space to
$\Def(X)$ is 
\[
H^1(X-\{0\},\T_X)=H^1(\tilde X-E,\T_{\tilde X})=H^2_Z(\T_X),
\]
by \cite[Thm.~2]{[42]}. The map 
$H^1(\tilde X,\T_{\tilde X})\rightarrow H^1(\tilde
X-E, \T_{\tilde X})$
is the differential of $\Def(\tilde X)\rightarrow \Def(X)$.
Hence the exact sequence
$$H^1(\tilde X,\T_{\tilde X})\rightarrow H^1(\tilde X-E,\T_{\tilde X})
\rightarrow H^2_E(\tilde X,\T_{\tilde X})\rightarrow H^2(\tilde X,\T_{\tilde X})$$
is identical to
$$H^0(R^1\pi_*\T_{\tilde X})\rightarrow H^2_Z(\T_X)\rightarrow H^2_E(\T_{\tilde X})
\rightarrow H^0(R^2\pi_*\T_{\tilde X})$$
which yields the desired sequence. Elements of
$\O_{X,0}$ pull back to elements of $\O_{\tilde X}$, which then act on $\T_{\tilde
X}$, and so
$H^1(\tilde X,\T_{\tilde X})$ and $H^1 (\tilde X-E,\T_{\tilde X})$ are naturally
$\O_{X,0}$-modules.
\end{proof}

\begin{say}
\label{(3.3)}
In our situation we will be interested in the case that $(X,0)$ is an isolated
rational Gorenstein point, and that $\tilde X\rightarrow X$ is a crepant
resolution. Recall from \cite{[38]} that there is an invariant $k$ associated with a
rational Gorenstein point as follows:
\begin{itemize}
\item[$k=0$] if $(X,0)$ is a cDV point, so is terminal.
\item[$k=1$] if $(X,0)$ is a hypersurface singularity locally of the form $x^2+y^3
+f(y,z,t)=0$ where $f=yf_1(z,t)+f_2(z,t)$ and $f_1$ (respectively $f_2$) is a sum
of monomials $z^at^b$ of degree $a+b\ge 4$ (respectively $\ge 6$). 
\item[$k=2$] if $(X,0)$ is a hypersurface singularity locally of the form
$x^2+f(y,z,t)=0$ where $f$ is a sum of monomials of degree $\ge 4$.
\item[$k\ge 3$] if $mult_0 X=k$ and $emb.dim.(X,0)=k+1$. The exceptional divisor
of the blow-up of $0\in X$ is a del Pezzo surface of degree $k$.
\end{itemize}

So in particular, for $k\le 3$, $(X,0)$ is a hypersurface singularity, and for
$k=4$, $(X,0)$ is a complete intersection (Gorenstein in codimension 2 implies
complete intersection) of two equations whose leading terms are quadratic and
define a del Pezzo surface in $\Pfour$.
However
$(X,0)$ is never a complete intersection if
$k>4$.
\end{say}

We next generalise the Bogomolov-Tian-Todorov Theorem to the case of a germ of a 
crepant resolution of a rational Gorenstein threefold singularity.

\begin{proposition}
\label{3.4} 
Let $(\tilde X,E) 
\rightarrow (X,0)$ be a crepant resolution of an
isolated rational Gorenstein threefold singularity $(X,0)$. Then $\Def(\tilde X)$
is non-singular.
\end{proposition}

\begin{proof}
The Hodge theory of $\tilde X$ is well-behaved above the middle dimension
from \cite{[31]}: in particular the spectral sequence $$H^q(\Omega^p_{\tilde
X})\Rightarrow H^n(\tilde X, \CC)$$ degenerates at the $E_1$ level for $p+q>3$.
Thus
as in the proof of \cite[5.5]{[4]}, if $\tilde X_n/A_n$ is a
deformation of
$\tilde X$ over $A_n$, then $H^2(\Omega^2_{\tilde X_n/A_n})$ is a locally free
$A_n$-module of rank $\dim H^2(\Omega^2_{\tilde X})$. Then, as a case of
Ran's ``$T^2$-injecting'' criterion \cite{[36]}, we see we have an exact sequence
$$H^1(\Omega^2_{\tilde X_n/A_n})\rightarrow H^1(\Omega^2_{\tilde X_{n-1}/A_{n-1}})
\rightarrow H^2(\Omega^2_{\tilde X})\mapright{\phi} H^2(\Omega^2_{\tilde X_n/A_n})
\mapright{} H^2(\Omega^2_{\tilde X_{n-1}/A_{n-1}})$$
with $\phi$ injective, so that $T^1(\tilde X_n/A_n)\rightarrow T^1(\tilde X_{n-1}/
A_{n-1})$ is always surjective. Thus, by the $T^1$-lifting criterion,
(Theorem \ref{1.8})  $\Def(\tilde
X)$ is smooth. 
\end{proof}

\begin{say}
\label{(3.5)}
Suppose furthermore that $X$ is a complete intersection singularity with
a crepant resolution $\pi:\tilde X\rightarrow X$, with
embedding dimension $e$ given by $f_1=\cdots=f_n=0$ with $f_i\in {\bf
C}\{x_1,\ldots,x_e\}$. So $\O_{X,0}=\CC\{x_1,\ldots,x_e\}/(f_1,\ldots,f_n)$.
(Of
course, given that
$X$ is a threefold canonical singularity, either $n=1,e=4$ or $n=2,e=5$.) Then
by \cite[pg.~634]{[44]},
$$T^1\cong \O_{X,0}^n/J$$
where $J$ is the submodule of $\O_{X,0}^n$ generated by $(\partial f_1/\partial x_i,
\ldots,\partial f_n/\partial x_i)$, $1\le i\le e$. The obvious $\O_{X,0}$-module
structure on $\O_{X,0}^n/J$ coincides with the $\O_{X,0}$-module structure
on $T^1$ from Lemma \ref{3.2}. If we
choose elements 
\[
(g_{11},\ldots,g_{1n}),\ldots,
(g_{m1},\ldots,g_{mn})\in\m\O_{X,0}^n
\]
which  along with
$(1,0,\ldots,0),\ldots,(0,\ldots,0,1)$ form a basis for $T^1$ after
reducing modulo $J$, then a miniversal family over the germ $\Def(X)=(T^1,0)$ about
the origin of
$T^1$ is given by
\begin{align*}
f_1+a_1+b_1g_{11}+\cdots+b_mg_{m1} ={} & 0\\
\vdots&\\
f_n+a_n+b_1g_{1n}+\cdots+b_mg_{mn} = {} & 0
\end{align*}
where $a_1,\ldots,a_n,b_1,\ldots,b_m$ are coordinates on $(T^1,0)$ given
by our choice of basis. This
defines a miniversal deformation $F:(\X,0)\rightarrow (T^1,0)$ of $(X,0)$. $F$ has a 
discriminant locus $D\subseteq (T^1,0)$, over which the fibres of $F$
are singular. From Lemma \ref{3.2}, we have a quotient $T'$ of $T^1$ which is an
$\O_{X,0}$-module.
\end{say}

A few comments about our plan to prove Theorem \ref{3.8} are in order here. The basic
strategy is to show that certain tangent vectors in $T^1$ always correspond to
smoothing directions. Lemma \ref{3.7} will show these tangent vectors will be the tangent
vectors which do not land in $\m T'$ under the projection $T^1\rightarrow T'$.
Lemma \ref{3.6} helps us identify tangent vectors for which this is not the case.

\begin{lemma}
\label{3.6} 
In the situation of \ref{(3.5)}, suppose $(X,0)$ is not an ordinary
double point (analytically isomorphic to $x_1^2+\cdots+x_4^2=0$).
Let ${\bf m}$ be the maximal ideal of $\O_{X,0}$. Then $T'\not=0$ and
if $x\in T'$ is annihilated by
${\bf m}$, then $x\in {\bf m}T'$.
\end{lemma}

\begin{proof}
We split this into the hypersurface case and the codimension two complete
intersection case, the former being simpler than the latter. First suppose
$(X,0)$ is a hypersurface singularity, so that $T'\cong\O_{X,0}/I$ for some ideal 
$I$ containing the jacobian ideal $J$. If some element of $T'$ not in ${\bf
m}T'$ were annihilated by ${\bf m}$, we would then have a non-zero element of
$T'/{\bf m}T'$ killed by ${\bf m}/{\bf m}^2$ under the surjective multiplication map
$T'/{\bf m}T'
\otimes {\bf m}/{\bf m}^2\rightarrow {\bf m}T'/{\bf m}^2T'$. But $\dim_{\CC}
T'/{\bf
m}T'
\le 1$, so this is only possible if ${\bf m}T'=0$ and $\dim_{\CC} T'=1$. Thus we
just need to show that $\dim_{\CC} T'>1$.

By Proposition \ref{3.4}, 
$\Def(\tilde X)$ is smooth. On the other hand, the fibre of the miniversal
space of $X$ over a general point of $D\subseteq \Def(X)$ has one ODP (\cite[\S5]{[44]}).
I claim that in fact the image of $\Def(\tilde X)$ in $\Def(X)$ cannot contain
any points corresponding to deformations of $X$ to a non-singular germ or a germ
with one ODP, and thus the codimension of the image of $\Def(\tilde X)$ in $\Def(X)$
is at least 2, from which we conclude $\dim T'\ge 2$. To show this claim, if
the invariant $k$ of $(X,0)$ is at least one, then $(\tilde X,0)$ contains some
exceptional divisors, and thus by \cite[Lem.~3.1]{[26]}, any deformation of $(\tilde
X,0)$ does also. Since any deformation of $(\tilde X,0)$ also has trivial canonical
bundle, these divisors blow down to yield a singularity which is not an ODP.
If $k=0$, $\tilde X\rightarrow X$ is a small resolution, and the claim follows from
\cite[Lem.~(1.8)]{[29]}.

In the case that $(X,0)$ is a codimension two complete intersection,
let 
\[
\tilde X
\mapright{\pi_1} X_1\mapright{\pi_2} X
\]
be a factorization of $\pi$, with $\pi_2$
the blowing up of $X$ at $0$.  By the argument of \cite[(1.5)]{[48]}, we have
$\ext^1(\Omega^1_{X_1},\O_{X_1})
\cong \ext^1(\pi_1^*\Omega^1_{X_1},\O_{\tilde X})$,
$\ext^1(\Omega^1_X,\O_X)
\cong \ext^1(\pi^*\Omega^1_X,\O_{\tilde X})$, and the natural maps $\pi^*\Omega^1_X
\rightarrow \pi_1^*\Omega^1_{X_1}\rightarrow\Omega^1_{\tilde X}$ induce the maps
on tangent spaces
$$\ext^1(\Omega^1_{\tilde X},\O_{\tilde X})\rightarrow
\ext^1(\pi_1^*\Omega^1_{X_1},\O_{\tilde X})\rightarrow 
\ext^1(\pi^*\Omega^1_X,\O_{\tilde X}).$$
Thus if 
$$T''=
\coker(\ext^1(\Omega^1_{X_1},\O_{X_1})\rightarrow \ext^1(\Omega^1_X,\O_X)),$$
there is a surjection $T'\rightarrow T''\rightarrow 0$ and it is enough to show
that $T'/\m T'\cong T''/\m T''$, $T''\not=0$,
and any element not in ${\bf m}T''$ is not
annihilated by
${\bf m}$. To do this, we consider $X\subseteq Y=({\CC}^5,0)$
and $Y_1\rightarrow Y$ the blowing-up
of the origin, $X_1\subseteq Y_1$ the proper transform of $X$ in $Y_1$. We will
first show that every infinitesimal deformation of $X_1$ is a deformation of $X_1$
inside of $Y_1$. We have an exact sequence
$$\exact{\I/\I^2}{\Omega^1_{Y_1}|_{X_1}}{\Omega^1_{X_1}}$$
with $\I$ the ideal sheaf of $X_1\subseteq Y_1$, which yields
$$\hom(\I/\I^2,\O_{X_1})\rightarrow\ext^1(\Omega^1_{X_1},\O_{X_1})
\rightarrow \ext^1(\Omega^1_{Y_1}|_{X_1},\O_{X_1}).$$
Since $\Omega^1_{Y_1}|_{X_1}$ is locally free, 
$\ext^1(\Omega^1_{Y_1}|_{X_1},\O_{X_1})=H^1(\T_{Y_1}|_{X_1})$, and we have an exact
sequence
$$H^1(\T_{Y_1})\rightarrow H^1(\T_{Y_1}|_{X_1})\rightarrow H^2(\T_{Y_1}\otimes\I).$$
Let $F$ be the exceptional $\Pfour$ of the blowing-up $Y_1\rightarrow Y$. 
Then using 
$$H^1(\T_{Y_1})=\lim_{\leftarrow} H^1(\T_{Y_1}\otimes\O_{Y_1}/\I_F^n),$$ the exact
sequences
$$\exact{\T_F}{\T_{Y_1}|_F}{\O_F(-1)}$$
and
$$\exact{\O_F(n-1)=\I_F^{n-1}/\I_F^n}{\O_{Y_1}/\I_F^n}{\O_{Y_1}/\I_F^{n-1}},$$
we see that $H^1(\T_{Y_1})=0$. To show that $H^2(\T_Y\otimes \I)=0$,
we use the Koszul resolution of $\I$
$$\exact{\L_1}{\L_2\oplus\L_3}{\I}$$
where $\L_1,\L_2$ and $\L_3$ are line bundles, with $\L_1|_F\cong\O_F(-4)$,
$\L_2|_F\cong\L_3|_F\cong\O_F(-2)$. We then see as before that $H^2((\L_2\oplus\L_3)
\otimes\T_{Y_1})=0$ and $H^3(\L_1\otimes\T_{Y_1})=0$.
Thus $H^2(\T_{Y_1}\otimes\I)=0$. To conclude, we have shown that
$\ext^1(\Omega^1_{Y_1}|_{X_1},\O_{X_1})=0$, and so any infinitesimal deformation of
$X_1$ comes from an  infinitesimal deformation of $X_1$ in $Y_1$. 

Now a choice of $g=(g_1,g_2)\in\O_{X,0}^2$ gives a deformation $\bar X/A_1$ of
$X$ by the equations
\begin{align*}
f_1+tg_1={} & 0\\
f_2+tg_2={} & 0.
\end{align*}
Note that if $\bar X/A_1$ is normally flat along $0\times \Spec A_1$ then
$(g_1,g_2)\in \m^2\O_{X,0}^2$. Indeed, since $f_1=f_2=0$ defines a rational
Gorenstein singularity, the leading terms of $f_1$ and $f_2$ must be quadratic,
and thus if $(g_1,g_2)\not\in\m^2\O_{X,0}^2$, the blow-up of $\bar X/A_1$ 
at $0\times\Spec A_1$ would not yield an exceptional divisor flat over $A_1$. 

If we write $T''=\O_{X,0}^2/I$, with $I\supseteq J$, then $(g_1,g_2)\in I$
if $(g_1,g_2)$ yields a deformation $\bar X/A_1$ which lifts to a deformation
$\bar X_1/A_1$ of $X_1$. Since $\bar X_1/A_1$ is equivalent to a deformation
$\bar X_1'/A_1\subseteq Y_1\times\Spec A_1$, we can blow down this deformation
to obtain a deformation $\bar X'/A_1$ equivalent to $\bar X/A_1$, which in
particular is normally flat along $0\times\Spec A_1$. Since two deformations of $X$
given by $(g_1,g_2)$ and $(g_1',g_2')$ are equivalent if and only if
they differ by an element
$(h_1,h_2)\in J$, we must have $(g_1+h_1,g_2+h_2)\in\m^2\O_{X,0}^2$ for some
$(h_1,h_2)\in J$. 

Now $J\subseteq\m\O_{X,0}^2$, so if $(g_1,g_2)\not\in\m\O_{X,0}^2$, we cannot have
$(g_1,g_2)
\in I$. Thus $\dim_{\CC}T''/\m T''=2$, and since $\dim_{\CC} T'/\m T'\le 2$,
we must have the surjection $T'/\m T'\rightarrow T''/\m T''$
being an isomorphism, and $T''\not=0$.

If $(g_1,g_2)\not\in\m\O_{X,0}^2$, but $x_i(g_1,g_2)\in I$ for all $i$, then in
particular for each $i$, there exists $(h_1,h_2)\in J$ such that
$x_i(g_1,g_2)+(h_1,h_2)\in\m^2\O_{X,0}^2$, from which one easily sees that the
leading (quadratic) terms of $f_1$ and $f_2$ have proportional partial derivatives,
and hence must be proportional themselves, ruling out the possibility of having a del
Pezzo surface as an exceptional divisor. Thus $(g_1,g_2)$ is not annihilated by all
elements of $\m$ in $T''$.
\end{proof}

\begin{lemma}
\label{3.7} 
Let $F':(\X',0)
\rightarrow(\SS',0)$ be a flat deformation of
$(X,0)$, with $\SS'$
non-singular at $0\in\SS'$ and the tangent space at $0\in\SS'$ the vector
space
$T$. Assume furthermore that $X$ is not an ordinary double point. Since
$F:(\X,0)\rightarrow (T^1,0)$ given in \ref{(3.5)} is a miniversal family for $(X,0)$,
there is a (non-unique) map 
$\SS'
\rightarrow (T^1,0)$ inducing a unique differential $T\rightarrow T^1$. Composing this
map with
$T^1\rightarrow T'$, we obtain a map $T\rightarrow T'$. If
$\im(T\rightarrow T')\not\subseteq {\bf m}T'$, then a general fibre of
$\X'\rightarrow\SS'$ is non-singular.
\end{lemma}

\begin{proof}
Following \cite[pg.~645]{[44]}, let $D\subseteq (T^1,0)$ be the
discriminant locus of $F$, $C\subseteq (\X,0)$ the critical locus. From \cite{[44]},
$F|_C:C\rightarrow D$ is the normalization of $D$. If every fibre of $F'$ is
singular, the induced map $\SS'\rightarrow (T^1,0)$ factors through $D$.
Since $\SS'$ is non-singular, in particular this map
factors through $C\rightarrow D$, and thus
the image of $T$ in $T^1$ is contained in the image in $T^1$ of the tangent space
$T_{\X,0}$ of
$(\X,0)$ at
$0$ via $F_*$. Using the explicit equations
for $(\X,0)$ and basis for $T^1$ given in \ref{(3.5)}, we see that the image of $T_{\X,0}$
in
$T^1$ is a vector space
$W\subseteq T^1$ given by
$a_1=\cdots=a_n=0$.
Again from the explicit description of $T^1$ in
\ref{(3.5)}, $W=\m T^1$, so the image of $W$ in $T'$ is a subspace
$V\subseteq \m T'$. Thus, we see that if $F'$ only
has singular fibres, the image of $T$ is contained in ${\bf m}T'$.
\end{proof}

\begin{theorem}
\label{3.8} 
Let $\tilde X$ be a non-singular Calabi-Yau threefold,
and $\pi:\tilde X\rightarrow X$ be a birational contraction morphism
such that $X$ has isolated, canonical, complete intersection singularities. Then
there is a deformation of $X$ which smooths all singular points of $X$
except possibly the ordinary double points of $X$. In particular, if $X$ has no
ordinary double points, then $X$ is smoothable.
\end{theorem}

\begin{proof}
Let $P\in Z=\Sing(X)$ be a singular point of $X$ which is not
an ordinary double point. We will show that there is a deformation of $X$
which smooths $P$. We have a diagram
\begin{equation}
\label{3.9}
\xymatrix@C=30pt
{
H^1(\tilde X,\T_{\tilde X})
\ar[r]& H^1(U,\T_{\tilde X}) \ar[r]&H^2_{\pi^{-1}(Z)}(\tilde X,
\T_{\tilde X})\ar[r]&H^2(\tilde X,\T_{\tilde X})\\
& H^1(U,\T_{X})\ar[u]^{\cong} \ar[r]&H^2_{Z}(X,
\T_{X})\ar[u]&
}
\end{equation}
$\Def(X)$ is smooth by \cite{[28]} or Theorem \ref{2.2}. The second vertical
map arises as in the proof of Lemma \ref{3.2}.
Thus to show that there is a
deformation of
$X$ which smooths
$P\in X$, it is enough to show by Lemma \ref{3.7} that the image of the
composed map 
$$H^1(U,\T_{\tilde X})\rightarrow H^2_{\pi^{-1}(Z)}(\tilde X,\T_{\tilde X})
\rightarrow H^2_E(\tilde X,\T_{\tilde X})$$ with $E=\pi^{-1}(P)$ contains
an element not in ${\bf m}_PT'_P
\subseteq T'_P\subseteq H^2_E(\tilde X,\T_{\tilde X})$. Here ${\bf m}_P$ is the 
maximal ideal of $\O_{X,P}$, and $T'_P$ is the subspace of $H^2_E(\tilde
X,\T_{\tilde X})$ given by Lemma \ref{3.2} applied to the germ $(X,P)$. To show this, it is
enough to show that there is an element of $\ker(H^2_E(\tilde X,\T_{\tilde X})
\mapright{\phi} H^2(\tilde X,\T_{\tilde X}))$ not in ${\bf m}_PT'_P$.

To see this, first consider the dual map $$H^1(\Omega^1_{\tilde X})
\mapright{\dual{\phi}} (R^1\pi_*\Omega^1_{\tilde X})_P.$$ Since 
$H^1(\Omega^1_{\tilde X})\cong (\Pic \tilde X)\otimes_{\ZZ}{\CC}$, this
map factors through $\Pic(\tilde X,E)\otimes_{\ZZ}{\CC}$, where $(\tilde X,E)$
denotes the germ of $\tilde X$ at $E$.
Now let
$E=\bigcup E_i$, with $E_i$ irreducible.

{\it Claim:} The map $\Pic(\tilde X,E)\otimes_{\ZZ} {\CC}
\rightarrow \bigoplus H^1(\Omega^1_{E_i})$ is injective.

\begin{proof}
Let $(\tilde X',E')\rightarrow (\tilde X,E)$ be a composition
of blow-ups of non-singular subvarieties so that $(\tilde X',E')
\rightarrow (X,P)$ is a resolution with $E'=\bigcup
E'_i$ simple normal crossings. Then we have a diagram
\[
\xymatrix@C=30pt
{
0\ar[r]&\Pic(\tilde X',E')\otimes_{\ZZ}\CC
\ar[r]&\bigoplus H^1(\Omega^1_{E_i'})\\
&\Pic(\tilde X,E)\otimes_{\ZZ}\CC\ar[u]
\ar[r]&\bigoplus H^1(\Omega^1_{E_i})\ar[u]\\
&0\ar[u]&
}
\]
The first column is exact since
$\tilde X'\rightarrow \tilde X$ is a series of blow-ups. The map
$\bigoplus H^1(\Omega^1_{E_i})\rightarrow\bigoplus H^1(\Omega^1_{E_i'})$
is defined via the induced maps $E_i'\rightarrow E_j$ for all $i,j$.
Suppose the
first row were exact. Then the claimed map must be injective.

The exactness of the first row follows from an argument of Namikawa in a preprint
version of \cite{[28]}; since this argument did not appear in the final version of the
paper, we sketch it here.

First, by \cite{[38]}, the $E_i'$ are all rational or ruled surfaces, since they are
exceptional divisors in the resolution of a rational Gorenstein point. Let
$\tau_{E'}^p\subseteq \Omega_{E'}^p$ be the torsion subsheaf. Then there is a
spectral sequence
$$E_1^{p,q}=H^q(E',\Omega_{E'}^p/\tau_{E'}^p)\Rightarrow H^{p+q}(E',{\CC})$$
which degenerates at the $E_1$ term, by \cite[Prop.~1.5]{[7]}. Note that
$H^{p+q}(E',{\CC})\cong
H^{p+q}((\tilde X',E'),{\CC})$, the cohomology of the germ $(\tilde X',E')$.
Let $$E'_{[p]}=\coprod_{i_0<\cdots<i_p} E_{i_0}'\cap\cdots\cap E_{i_p}'.$$
There is an exact sequence
$$0\rightarrow \Omega_{E'}^p/\tau_{E'}^p\rightarrow
\Omega^p_{E'_{[0]}}\rightarrow\Omega^p_{E'_{[1]}}\rightarrow\cdots$$ for each $p$,
again by \cite[Prop.~1.5]{[7]}.
First consider this sequence for $p=0$. Since $H^2(\O_{E'})=0$, as $(X,P)$ is
a rational singularity, we find $H^1(\O_{E'_{[0]}})\rightarrow H^1(\O_{E'_{[1]}})$
is surjective, and thus $H^0(\Omega^1_{E'_{[0]}})\rightarrow
H^0(\Omega^1_{E'_{[1]}})$ is surjective.

By the same sequence for $p=1$, we find $H^1(\Omega_{E'}^1/\tau_{E'}^1)
\rightarrow H^1(\Omega^1_{E'_{[0]}})$ is injective. With $p=2$, we obtain
$H^0(\Omega_{E'}^2/\tau_{E'}^2)=H^0(\Omega^2_{E'_{[0]}})=0$ since all components
of $E'$ are ruled. Thus our spectral sequence yields $H^2(E',{\CC})\cong
H^1(\Omega_{E'}^1/\tau_{E'}^1)$ and so the map $H^2(E',{\CC})\rightarrow
H^1(\Omega^1_{E'_{[0]}})$ is injective. Since $\Pic(\tilde X',E')\otimes_{\ZZ}
{\CC}\cong H^2((\tilde X',E'),{\CC})\cong H^2(E',{\CC})$, the result follows.
\end{proof}

Thus the composed map 
$$\Pic(\tilde X,E)\otimes_{\ZZ}{\CC}\rightarrow
(R^1\pi_*\Omega^1_{\tilde X})_P\rightarrow \bigoplus_i  H^1(\Omega^1_{\tilde
X}|_{E_i})$$ must also be injective, as is then the composed map
$$\coim(\dual{\phi})\rightarrow (R^1\pi_*\Omega^1_{\tilde X})_P
\rightarrow\bigoplus_i H^1(\Omega^1_{\tilde X}|_{E_i}).$$
Dually, we get that the composed map 
$$\bigoplus_i \ext^2_{\O_{\tilde X}}(\O_{E_i},\T_{\tilde X})
\rightarrow H^2_E(\T_{\tilde X})\rightarrow \im(\phi)$$
is surjective. Now clearly $\ext^2_{\O_{\tilde X}}(\O_{E_i},\T_{\tilde X})$
is annihilated by the maximal ideal
${\bf m}_P\subseteq\O_{X,P}$ since $\O_{E_i}$ is. Let
$W\subseteq
\bigoplus \ext^2_{\O_{\tilde X}}(\O_{E_i},\T_{\tilde X})$ be a subspace mapping 
isomorphically via $\phi$ to $\im(\phi)$. Identifying $W$ with its image
in $H^2_E(\T_{\tilde X})$, we must then have $H^2_E(\T_{\tilde X})
=W+\ker(\phi)$. $W$ is annihilated by ${\bf m}_P$, and so $W\cap T'_P\subseteq
{\bf m}_PT'_P$ by Lemma \ref{3.6}; since $T_P'\not=0$,  $\ker(\phi)$ contains some elements not in ${\bf
m}_PT'_P$.
\end{proof}

Combining this with Namikawa and Steenbrink's results in \cite{[30]}, we obtain

\begin{corollary}
\label{3.10} Let $\tilde X$ be a $\QQ$-factorial Calabi-Yau threefold with
terminal singularities, and suppose $\pi:\tilde X\rightarrow X$ is a birational
contraction morphism such that $X$ has isolated, canonical, complete intersection
singularities. Then there is a deformation of $X$ to a variety with at worst
ordinary double points.
\end{corollary}

\begin{proof}
 By \cite{[30]}, there is a small deformation of $\tilde X$, $\tilde\X\rightarrow\Delta$,
such that $\tilde\X_t$, $t\not=0$ is non-singular.
If $H$ is an ample Cartier divisor on $X$,
$\pi^* H$ is a nef and big divisor on $\tilde X$, and by \cite[Prop.~6.1, Thm.~C]{[28]},
this divisor deforms to a relatively nef and big divisor on $\tilde \X$. Thus
the morphism $\pi$ deforms to a morphism
$\pi':\tilde \X\rightarrow \X$, with $\X$ a deformation of $X$. (It is possible
this may introduce some new ordinary double points on $\X_t$.)
$\X_t$ still has isolated
complete intersection canonical singularities, but has
a crepant resolution $\tilde\X_t\rightarrow\X_t$, so we can apply Theorem \ref{3.8}.
\end{proof}

\section{Calabi-Yaus with non-complete intersection singularities}
\label{sec4}

If $X$ is a Calabi-Yau with non-complete intersection isolated canonical
singularities, then there is no statement as complete as Theorem \ref{3.8}.
There are three difficulties. First, not every isolated singularity is smoothable.
Second, as seen in Example \ref{2.7}, the deformation theory of the Calabi-Yau can be quite
bad. Third, as Example \ref{4.1} shows, even if the singularity is smoothable, we may not
even have any infinitesimal deformations of $X$ which yield a non-trivial
infinitesimal deformation of the singularity.

\begin{example}
\label{4.1} Let $S$ be a non-singular del Pezzo surface of degree $\ge 5$,
and consider an elliptic fibration $f:\tilde X\rightarrow S$ defined by the 
Weierstrass equation $y^2=x^3+ax+b$, for general $a\in H^0(\omega_S^{-4})$,
$b\in H^0(\omega_S^{-6})$. $\tilde X$ will be a non-singular Calabi-Yau threefold,
and $f$ has a section $\sigma$, the section at infinity obtained after
compactifying the affine equation given. Identifying $S$ with $\sigma(S)$, it is
possible to find a contraction $\pi:\tilde X\rightarrow X$ contracting $S$. As in
\eqref{3.9}, we
have the exact sequence
$$H^1(\tilde X,\T_{\tilde X})\rightarrow H^1(U,\T_{\tilde X})
\rightarrow H^2_S(\tilde X,\T_{\tilde X}) \rightarrow H^2(\tilde X,\T_{\tilde X})$$
where $U=\tilde X-S$. Dualizing the last map we obtain the map
$H^1(\tilde X,\Omega^1_{\tilde X})\rightarrow H^0(R^1\pi_*\Omega^1_{\tilde X})$.
Using $\deg S\ge 5$, a straightforward calculation shows that
\[
H^0(R^1\pi_*\Omega^1_{\tilde X})\cong H^1(S,\Omega^1_S),
\]
and $H^1(\tilde X,
\Omega^1_X)\rightarrow H^1(S,\Omega^1_S)$ is the restriction map $\sigma^*$ and is
thus surjective, since $\sigma^*f^*$ is the identity on $H^1(S,\Omega^1_S)$. Thus
$T^1(\tilde X/k)=H^1(\tilde X,\T_{\tilde X})
\rightarrow T^1(X/k)=H^1(U,\T_{\tilde X})$ is surjective, so
there are no infinitesimal deformations of $X$ which don't 
come from infinitesimal deformations of $\tilde X$. Thus 
$\Def(\tilde X)\rightarrow \Def(X)$ is surjective.
Since the exceptional locus $S$ deforms in any deformation of $\tilde X$,
$X$ is not smoothable.
\end{example}

It is clear from this example that one issue is controlling the map $T^1\rightarrow
T^1_{\loc}$, where $T^1$ is the tangent space to $\Def(X)$ and $T^1_{\loc}$ is the
tangent space to $\Def(X,Z)$, $Z=\Sing(X)$. An assumption which will give us as much
control over this map as possible is that $X$ is $\QQ$-factorial. Nevertheless,
this does not guarantee surjectivity of $T^1\rightarrow T^1_{\loc}$. Hence I make
here some further assumptions which I feel are quite artificial and which hopefully
can be removed in the future, once more is known about the deformation theory of
rational Gorenstein threefold singularities.

We make the following rather ad hoc definition:

\begin{definition}
\label{4.2} An isolated non-complete intersection rational Gorenstein
point
$(X,P)$ is {\it good} when
\begin{enumerate}
\item If $X'\rightarrow X$ is the blow-up of $X$ at $P$ with exceptional
divisor $E$, then $E$ is irreducible and $X'$ has only isolated singularities.
\item $\Def(X')$ is non-singular, and the natural map $\Def(X')\rightarrow
\Def(X)$ is an immersion. 
\item There is a smoothing component of $\Def(X)$ containing the image of
$\Def(X')\rightarrow \Def(X)$. 
\end{enumerate}
\end{definition}

This is a particularly strong set of assumptions. In \S 5, we will prove certain
singularities are good. Here, we just want to make explicit the assumptions we need.

\begin{theorem} 
\label{4.3} 
Let $\tilde X$ be a non-singular Calabi-Yau threefold and
$\pi:\tilde X\rightarrow X$ a birational contraction such that $X$ is 
$\QQ$-factorial and for each $P\in \Sing(X)$, $(X,P)$ is good.
Then $X$ is smoothable.
\end{theorem}

First we need

\begin{lemma} 
\label{4.4} 
If $X$ is a compact, $\QQ$-factorial algebraic 
variety with rational singularities, $\pi:\tilde X\rightarrow X$
a resolution of singularities with irreducible $E_1,\ldots,E_n\subseteq \tilde X$ the
$\pi$-exceptional divisors, and $H^2(\O_{\tilde X})=0$, then
$$\im\left[ H^1(\Omega^1_{\tilde X})\mapright{p_1} H^0(R^1\pi_*\Omega^1_{\tilde X})
\right]
=\im\left[\bigoplus {\CC}E_i\mapright{p_2} H^0(R^1\pi_*\Omega^1_{\tilde X})
\right]$$
with the map $p_2$ the composition of $\bigoplus{\CC}E_i\rightarrow
H^1(\Omega^1_{\tilde X})$ and $p_1$. 
\end{lemma}

\begin{proof}
\cite[12.1.6]{[23]} gives a similar statement for rational cohomology,
i.e. $X$ $\QQ$-factorial implies that 
$$\im\left[ H^2(\tilde X,\QQ)\mapright{p_1} H^0(X,R^2\pi_*\QQ)
\right]
=\im\left[\bigoplus \QQ E_i\mapright{p_2} H^0(X,R^2\pi_*\QQ)
\right].$$
We can clearly replace $\QQ$ by ${\CC}$, and since $H^2(\O_{\tilde X})=0$,
$H^2(\tilde X,{\CC})=H^1(\tilde X,\Omega^1_{\tilde X})$. Now 
$R^i\pi_*\O_{\tilde X}=0$ for $i>0$ since $X$ has rational
singularities. Thus by the spectral sequence $$R^q\pi_*\Omega^p_{\tilde X}\Rightarrow
R^n\pi_*{\CC},$$ there is a natural map
$$H^0(X,R^2\pi_*{\CC})\rightarrow H^0(X,R^1\pi_*\Omega^1_{\tilde X}),$$
yielding a commutative diagram
\[
\xymatrix@C=30pt
{
H^2(\tilde X,{\CC})\ar[r]^{\cong}\ar[d]&H^1(\tilde X,\Omega^1_{\tilde X})
\ar[d]_{p_1}\\
H^0(X,R^2\pi_*{\CC})\ar[r]&H^0(X,R^1\pi_*\Omega^1_{\tilde X})
}
\]
This gives the desired result. 
\end{proof}

The need for the $\QQ$-factorial hypothesis comes from the following Lemma. In general
it is difficult to show that the map $T^1\rightarrow T^1_{\loc}$ is surjective (see Remark \ref{4.7}),
but we can make do with

\begin{lemma} \label{4.5} 
Let $\tilde X$ be a non-singular Calabi-Yau threefold, 
$\pi:\tilde X\rightarrow
X$ a birational contraction such that $X$ is $\QQ$-factorial with isolated 
singularities. Let $Z=\Sing(X)$, $E=\pi^{-1}(Z)$. Let $T^1$, $T^1_{\loc}$ be as in \ref{(2.1)},
the tangent spaces of $\Def(X)$ and $\Def(X,Z)$ respectively, and let
$$T'=\ker(H^2_E(\T_{\tilde X})\rightarrow H^0(R^2\pi_*\T_{\tilde X}))$$
as in Lemma \ref{3.2}. Then the composed map
$$T^1\rightarrow T^1_{\loc}\rightarrow T'$$
is surjective.
\end{lemma}

\begin{proof}
First we claim that
$$\ker(H^2_E(\T_{\tilde X})\rightarrow H^0(R^2\pi_*\T_{\tilde X}))
=\ker(H^2_E(\T_{\tilde X})\rightarrow H^2(\T_{\tilde X})).$$
That the first space contains the second is clear, since the first map factors
through the second. Dualizing this statement, we need to prove that
$$\coker(H^1_E(\Omega^1_{\tilde X})\mapright{p_1} H^0(R^1\pi_*\Omega^1_{\tilde X}))
=\coker(H^1(\Omega^1_{\tilde X})\mapright{p_2} H^0(R^1\pi_*\Omega^1_{\tilde X})).$$
The first surjects on the second, i.e. $\im(p_1)\subseteq
\im(p_2)$. To show equality, we need to show $\im(p_1)\supseteq\im(p_2)$. Let
$\{E_i\}$ be the irreducible components of $E$. We have the diagram
\[
\xymatrix@C=30pt
{
&&0\ar[d]&&\\
&\bigoplus{\CC}E_i\ar[r]\ar[d]&\Pic\tilde X\otimes_{\ZZ}{\CC}\ar[r]
\ar[d]_{d\log_X}&\Pic U\otimes_{\ZZ}{\CC}\ar[r]\ar[d]_{d\log_U}&0\\
0\ar[r]&H^1_E(\Omega^1_{\tilde X})\ar[r]&H^1(\Omega^1_{\tilde X})\ar[r]\ar[d]&
H^1(U,\Omega^1_{\tilde X})&\\
&&0&&
}
\]
Here, $H^0(U,
\Omega^1_{\tilde X})=H^0(U,\Omega^1_X)=0$ by \cite{[43]}. This shows
that $\im(\bigoplus {\CC}E_i\rightarrow H^1(\Omega^1_{\tilde X}))
\subseteq \im(H^1_E(\Omega^1_{\tilde X})\rightarrow H^1(\Omega^1_{\tilde X}))$, and
thus by Lemma \ref{4.4} and the hypothesis that $X$ is $\QQ$-factorial, $\im(p_2)
\subseteq \im(p_1)$. 

Now diagram \eqref{3.9} is still valid in this situation, and we have by definition
$$T'=\ker(H^2_E(\T_{\tilde X})\rightarrow H^0(R^2\pi_*\T_{\tilde X})),$$
so identifying $T^1$ with $H^1(U,\T_X)$, \eqref{3.9} shows that 
$$\im(T^1\rightarrow T')=\ker(H^2_E(\T_{\tilde X})\rightarrow H^2(\T_{\tilde X})).$$
The result then follows. 
\end{proof}

The reason for condition (1) in the definition of good singularity is the following
lemma, which shows that we can ``smooth in one step.''

\begin{lemma}
\label{4.6} 
Let $\X\rightarrow\Delta$ be a one-parameter deformation of
a good singularity $(X,P)$. This induces a map $\Delta\rightarrow \Def(X)$. Then
either $\im(\Delta\rightarrow \Def(X))\subseteq \im(\Def(X')\rightarrow \Def(X))$ or
$\X_t$ has only hypersurface singularities for general $t\in\Delta$.
\end{lemma}

\begin{proof}
Suppose that for general $t\in\Delta$, $\X_t$ has worse than hypersurface
singularities. $\X$ must be singular at each point where $\X_t$ has worse than
hypersurface singularities, so we can assume there is a curve $C\subseteq \X$
dominating $\Delta$ along which $\X$ is singular. By making a suitable base change,
we can assume $C$ is the image of a section $\Delta\rightarrow\X$. It is enough to
show that $\X$ is normally flat along $C$; this will permit us to blow up
$\X$ along $C$ to obtain a deformation $\X'\rightarrow \Delta$ of $X'$. To show
that $\X$ is normally flat along $C$, it is enough to show that the multiplicity
of $\X_t$ at $C_t$ is constant for $t\in\Delta$. Indeed, $\X$ is normally flat
along $C$ if the Hilbert-Samuel function of the point $C_t\in\X_t$ is constant for
$t\in\Delta$. For a rational Gorenstein non-cDV point, the exceptional locus
upon blowing-up is a del Pezzo surface, and the Hilbert-Samuel function for a
family of del Pezzo surfaces of the same degree is constant.
 To show that the
multiplicity of
$\X$ along $C$ is constant, it is enough to show the same thing for a general
hyperplane section of $\X$ containing $C$. 

To this end, take a general hyperplane section of $\X$ containing $C$; this yields
a deformation $\SS\rightarrow\Delta$ of surface singularities, with $\SS_0$ a general
hyperplane section of $\X_0=X$. $\SS_0$ is then a Gorenstein elliptic singularity,
and
the blowing up of $\SS_0$ at $P$ resolves the singularity, since $(X,P)$ is good.
In particular, if $\tilde\SS_0\rightarrow\SS_0$ is the blowing up of $\SS_0$, then 
the exceptional locus is either a non-singular elliptic curve, a nodal rational
curve, or a cuspidal rational curve. These are the singularities $El(d)$, 
$C(d)$ or $Cu(d)$, where $d$ is the multiplicity of the singularity,
following the notation of \cite[\S5]{[49]}. $\SS_t$ must also have elliptic singularities.
If the list of adjacencies given in \cite[\S5]{[49]}, then $\SS_0$ can only
be adjacent to an elliptic singularity of the same multiplicity (these are
the adjacencies $Cu(d)\rightarrow C(d)$ and $C(d)\rightarrow El(d)$ in Wahl's notation).
In the case at hand, Wahl's list of adjacencies has been proved to be complete
in \cite{[27]}. Thus the multiplicity of the singularities of $\SS_t$ remain constant,
proving the lemma.
\end{proof}

\begin{proof}[Proof of Theorem \ref{4.3}] 
We will now prove the theorem by constructing a deformation $\X\rightarrow\Delta$
of $X$ such that $\im(\Delta\rightarrow \Def(X,Z))$ is contained in the smoothing
component given by item (3) of Definition \ref{4.2}, but
$$\im(\Delta\rightarrow \Def(X,Z))\not\subseteq \im(\Def(X',E')\rightarrow \Def(X,Z)),$$
where $X'\rightarrow X$ is the blow-up of $X$ at the singular locus $Z$ of $X$,
with exceptional locus $E'$. It will then follow from Lemma \ref{4.6} that
$Y=\X_t$ has only hypersurface singularities for $t\not=0$.
By \cite[12.1.11]{[23]}, $Y$ is $\QQ$-factorial. By Corollary \ref{3.10}, 
we can deform $Y$ to something with only ordinary
double points which is still $\QQ$-factorial, and hence can be smoothed
entirely by \cite{[8]}.

To show this, we will apply Theorems \ref{2.2} and \ref{1.9} to get dimension estimates
of components of $\Def(X)$ versus $Def(X')$. We adopt the following
notation:
\begin{align*}
S:={} & \hbox{the ring which pro-represents the deformation functor of $X$.}\\
\Lambda:={} & \hbox{a hull of the deformation functor of $(X,Z)$.}\\
\Lambda':={} & \hbox{a hull of the deformation functor of $(X',E')$.}\\ 
\tilde \Lambda:= {} & \hbox{a hull of the deformation functor of $(\tilde X,E)$.}\\ 
V_0:= {} & \ker(T^1\rightarrow T^1_{\loc}).\\ 
V_1:={} & \coker(T^1\rightarrow T^1_{\loc}).\\ 
r:={} & \dim V_0.\\ 
s:={} & \dim V_1.\\ 
R:={} & \Lambda[[x_1,\ldots,x_r]].
\end{align*}
By Theorems \ref{2.2} and \ref{1.9},
we have $S\cong R/J$ for some ideal $J\subseteq R$.
Furthermore, there is
an ideal $J'\subseteq J$ with $\Supp(R/J)=\Supp(R/J')$, $\m_R/(\m_R^2+J)\cong
\m_R/(\m_R^2+J')$ and $J'$ is generated by
$s$ elements.

In addition, we set 
\begin{align*}
p:\Spec R\rightarrow \Spec\Lambda&\quad\hbox{the projection.}\\
\scrC\subseteq \Spec R&\quad\hbox{the pull-back via $p$ of the smoothing component
given by item (3)}\\
&\quad\hbox{ of Definition \ref{4.2}.}\\
\scrB\subseteq \Spec R&\quad\hbox{the subscheme defined by $J'$, set-theoretically
equal to $\Def(X)$.}\\
\tilde\D&:=\Spec\tilde\Lambda\times_{\Spec\Lambda} \Spec R.\\ 
f:\tilde\D\rightarrow\Spec R&\quad\hbox{the projection.}\\ 
\D'&:=\Spec\Lambda'\times_{\Spec\Lambda} \Spec R.
\end{align*}
By item
(2) of Definition \ref{4.2}, $\D'$ can be identified with a non-singular subscheme
of $\Spec R$, and $\im f\subseteq\D'$. 

Finally, let $T_{\B,0}$, $T_{\D',0}$ and $T_{\tilde\D,0}$ be the Zariski
tangent
spaces to $\B$, $\D'$ and $\tilde\D$ at their closed points. $T_{\B,0}$ and $T_{\D',0}$
are contained in
$T=T_{\Spec R,0}=T^1_{\loc}\oplus V_0$. 
Then $f_*(
T_{\tilde\D,0})\subseteq T_{\D',0}$, and by Lemma \ref{3.2} $f_*(T_{\tilde \D,0})$ is the
kernel of the composed map $T\rightarrow T^1_{\loc}\rightarrow T'$.

Note that since $\m_R/(\m_R^2+J)\cong
\m_R/(\m_R^2+J')$, $T_{\B,0}\cong T^1$. Thus
\begin{align*}
\dim (T_{\B,0}+T_{\D',0})\ge {} &  \dim(T_{\B,0}+f_*(T_{\tilde\D,0}))\\
= {} & \dim T_{\B,0}
+\dim f_*(T_{\tilde\D,0})-\dim T_{\B,0}\cap f_*(T_{\tilde\D,0})\\
= {} & \dim T'+
\dim
f_*(T_{\tilde\D,0}),
\end{align*}
since
$f_*(T_{\tilde\D,0})\subseteq T_{\D',0}$ and Lemma \ref{4.5} gives a surjection
$T^1\rightarrow T'$ with kernel $f_*(T_{\tilde\D,0})\cap T^1$. 
Also,
$s=\dim T-\dim T_{\B,0}=\dim f_*(T_{\tilde\D,0})+\dim T'-\dim T_{\B,0}$. Then
\begin{align*}
\dim\B\cap\D'\le {} & \dim T_{\B,0}\cap T_{\D',0}\\
= {} & \dim T_{\B,0}+\dim T_{\D',0}-\dim(T_{\B,0}+T_{\D',0})\\
\le {} & \dim T_{\B,0}+\dim T_{\D',0}-\dim f_*(T_{\tilde\D,0})-\dim T'\\
= {} & \dim T_{\D',0}-s\\
={} & \dim\D'-s.
\end{align*}
On the other hand,
\begin{align*}
\dim\scrB\cap\scrC\ge {} & \dim\scrC-\hbox{\# of equations generating $J'$}\\
={} & \dim \scrC-s\\
\ge {} & \dim \D'+1-s\quad\quad\hbox{(since $\D'$ is strictly contained in $\scrC$)}\\
> {} &\dim\B\cap\D'.
\end{align*}
Thus $\B\cap\D'$ is strictly contained in $\scrB\cap\scrC$. Thus there is a map $\Spec k[[t]]
\rightarrow \B\cap\scrC$ whose image is not contained in $\B\cap \D'$. This yields
a small deformation $\X\rightarrow\Delta$ which has the desired properties.
\end{proof}

\begin{remark}
\label{4.7} 
More generally, if every singularity of $X$ is smoothable and,
in the notation of \ref{(2.1)}, $T^1\rightarrow T^1_{\loc}$ is surjective, then it is clear
that $X$ is smoothable, since then by Theorems \ref{2.2} and \ref{1.9}, any deformation
of $(X,Z)$ is realisable via a deformation of $X$. This will hold, for example, if the map
$H^0(R^1\pi_*\T_{\tilde X}) \rightarrow H^2_Z(\T_X)$ of Lemma \ref{3.2} 
is zero and $X$
is $\QQ$-factorial. 
\end{remark}

\section{Primitive Contractions}
\label{sec5}

Recall from \cite{[50]} that if $\tilde X$ is a non-singular Calabi-Yau, then
$\pi:\tilde X\rightarrow X$ is a {\it primitive contraction} if $\pi$ cannot be
factored in the algebraic category. The goal of this section is to apply the
results of the previous sections to the case that $\pi$ is a primitive contraction.
We will find strong restrictions on the possible primitive contractions of $\tilde
X$ if we assume that
$\tilde X$ is primitive in the sense of the introduction. Recall from 
\cite{[50]} the
following classification of primitive contractions:
\begin{itemize}
\item[Type I:] $\pi$ contracts a union of curves.
\item[Type II:] $\pi$ contracts a divisor to a point.
\item[Type III:]
$\pi$ contracts a divisor to a curve.
\end{itemize}

We treat the first two cases in this section. The type III case will be treated in
\cite{[12]}. Type I contractions have already been treated by Namikawa:

\begin{theorem}
\label{5.1} 
Suppose $\pi:\tilde X\rightarrow X$ is a primitive type I
contraction. Then $X$ is smoothable unless $\pi$ is the contraction of a single
$\Pone$ to an ordinary double point.
\end{theorem}

\begin{proof}
 If $X$ has only ordinary double points, and $C_1,\ldots,C_n$ are the
exceptional curves of $\pi$, then the cohomology classes $[C_1],\ldots,[C_n]\in
H_2(\tilde X,\ZZ)$ coincide, since $\pi$ is primitive, and thus unless $n=1$,
there
is a non-trivial linear dependence relation on $[C_1],\ldots, [C_n]$. Thus by
\cite{[8]}, $X$ is smoothable. If $X$ does not have only ordinary double points,
then let $Z\rightarrow X$ be a (non-projective) small resolution of the ordinary
double points of $X$. Then if $C_1,\ldots,C_n\subseteq Z$ are the exceptional
curves, $[C_1],\ldots,[C_n]=0$ in $H_2(Z,\ZZ)$, and so by \cite[Thm.~2.5]{[29]},
$X$ is smoothable.
\end{proof}

For type II contractions, we will need a more refined classification.

\begin{theorem}
\label{5.2} 
Let $\pi:\tilde X\rightarrow X$ be a primitive type II
contraction. Then if $X$ has a non-hypersurface singularity, the exceptional
divisor $E$ of $\pi$ is either
\begin{itemize}
\item[(i)] a normal rational del Pezzo surface of degree $\le 9$ or
\item[(ii)] a nonnormal del Pezzo surface of degree 7.
\end{itemize}
\end{theorem}

Proof: By \cite[Thm.~2.11]{[38]}, $\pi$ is the blowing-up of $X$ at $P$, the singular
point of $X$. The exceptional surface $E$ is a generalized del Pezzo surface
($\omega_E$ is ample) of degree
$k$, where $k$ is Reid's invariant (see \S 3). Since $\pi$ is primitive, $E$ is
integral.

We then have the following possibilities for $E$:
\begin{enumerate}
\item[(i)] $E$ is a normal, rational del Pezzo surface, in which case $\deg E\le
9$.
\item[(ii)] $E$ is a nonnormal del Pezzo surface, as classified in \cite{[39]}. The
possibilities are:
\begin{enumerate}
\item[(a)] Let $F_a=\P(\O_{\Pone}\oplus\O_{\Pone}(-a))$ be the rational scroll,
with $\Pic F_a$ generated by $C_0$, the negative section, and $f$, the class of a
fibre. Embed $F_a$ in $\P^{a+5}$ via $|C_0+(a+2)f|$. $E$ is the projection of $F_a$
into $\P^{a+4}$ from a point in a plane spanned by the conic $C_0$. This
projection maps $C_0$ two-to-one to a line $l$, and makes no other identifications.
\item[(b)] Embed $F_a\subseteq \P^{a+3}$ via $|C_0+(a+1)f|$. $E$ is the
projection of $F_a$ into $\P^{a+2}$ from a point in the plane spanned by the line
$C_0$ and one fibre $f$. This projection identifies $C_0$ and $f$.
\item[(c)] Take $E$ to be a cone over a rational nodal or cuspidal curve of
degree $d$ spanning $\P^{d-1}$. For $d>3$, the vertex of this cone will not be a
hypersurface singularity, and hence such a del Pezzo surface cannot be contained in
a non-singular Calabi-Yau threefold. Thus this case does not occur.
\end{enumerate}
\item[(iii)] $E$ is a cone over an elliptic normal curve, in which case this does
not occur just as in the cones in case (ii) (c).
\end{enumerate}
 
Thus we need to deal with cases (ii), (a) and (b). Suppose $D\subseteq X$ is such a
del Pezzo surface, $X$ a non-singular Calabi-Yau. Let $\tilde D\rightarrow D$ be
the normalization; $\tilde D$ is a scroll. Let $i:\tilde D\rightarrow X$ be the
induced map. Let $S\subseteq \tilde D$ be the subscheme defined by the zeroth
Fitting ideal of $\Omega^1_{\tilde D/X}$. This is the subscheme defined by the
condition that $i^*\Omega^1_X\mapright{di}\Omega^1_{\tilde D}$ drops rank, and so
is defined by the $2\times 2$ minors of the map $di$. (See \cite[III A]{[21]}.)
The degree of $S$ is the number of pinch points of $D$. Now in
case (ii) (a),
$S$ consists precisely of
the two pinch points corresponding to the ramification
points of
$C_0\rightarrow l$, and $S$ is a length two scheme. In case (ii) (b), there is just
one similar such point, but it is not an ordinary pinch point. To
analyze this point, we can consider $\AA^2\subseteq\tilde D$ with coordinates
$u$ and $v$ in such a way that $C_0$ and $f$ coincide with $u=0$ and $v=0$
respectively. The map $\tilde D\rightarrow D$ then identifies the $u$ and $v$ axes.
Following the recipe of \cite[2.1]{[39]}, $\tilde D\rightarrow D$ then locally looks like
$$(u,v)\in\AA^2\mapsto (u+v,uv,uv^2)\in
\AA^3.$$  The Jacobian
of this map is
\[
\begin{pmatrix}
1&v&v^2\cr 1&u&2uv
\end{pmatrix}
\]
and the ideal of $2\times 2$ minors is
$$(u-v,2uv-v^2,uv^2)=(u-v,uv),$$
which defines a scheme of length two at the origin. Thus in either case, $\deg S=2$.

Alternatively, we can compute $\deg S$ by \cite[(III, 8)]{[21]},
$$\deg S=c_2(i^*\T_X-\T_{\tilde D}).$$
Now the Chern polynomial in $t$ of $\T_{\tilde D}$ is
$$c_t(\T_{\tilde D})=1-K_{\tilde D}t+4t^2$$
and of $i^*\T_X$
$$c_t(i^*\T_X)=1+c_2(X).Dt^2.$$
We compute $c_2(X).D$ using Riemann-Roch:
$$\chi(\O_X(D))={1\over 6}D^3+{1\over 12}c_2(X).D.$$
Now $H^i(\O_X(D))=H^{3-i}(\O_X(-D))$ by Serre duality, and the exact sequence
$$\exact{\O_X(-D)}{\O_X}{\O_D}$$ shows that $H^i(\O_X(D))=0$ for $i>0$. 
($H^i(\O_X)=0$ for $i=1,2$ and $H^1(\O_D)=0$ by \cite[4.10]{[39]}.) Also,
$H^0(\O_X(D))=1$, so $\chi(\O_X(D))=1$. Thus 
$$c_2(X).D=12-2D^3.$$
We then find that
$$c_t(i^*\T_X-\T_{\tilde D})=c_t(i^*\T_X)/c_t(\T_{\tilde D})=1+K_{\tilde D}t
+(c_2(X).D+4)t^2$$
so $\deg S=16-2D^3$. Thus we must have $2=16-2D^3$, or $D^3=7$. Thus $D$ must be a
del Pezzo surface of degree 7. $\bullet$

\begin{remark}
\label{5.3}
While the above theorem was stated for the non-hypersurface
singularity case, in the hypersurface singularity case
we can still rule out the
non-normal case as above, but cones over elliptic curves can, and do, occur.
\end{remark}

\begin{proposition}
\label{5.4} 
Suppose $(X,0)$ is an isolated rational Gorenstein
threefold point with $k=mult_0 X\ge 5$, such that if $\tilde X\rightarrow X$ is the
blowing-up of
$X$ at $0$, then $\tilde X$ is non-singular and the exceptional divisor $E$ is
non-singular. Then $(X,0)$ is analytically isomorphic to a cone over $E$.
\end{proposition}

\begin{proof}
Let $\I_E$ be the ideal sheaf of $E$ in $\tilde X$. If $H^1(\T_E\otimes
(\I_E/\I_E^2)^n)=H^1((\I_E/\I_E^2)^n)=0$ for all
$n\ge 1$, then \cite[Cor.~to Satz 7]{[9]}, tells us that $(\tilde X,E)$ is
analytically isomorphic to an open neighborhood of $E$ embedded in the normal
bundle of $E$ in $\tilde X$ as the zero section. This will then give the theorem.

$E$ is a del Pezzo surface of degree between 5 and 9, and
$\I_E/\I_E^2=\omega^{-1}_E$. For a del Pezzo surface, $H^1(\omega_E^{-n})=0$ for
all $n\ge 1$. If we show that $H^1(\T_E)=0$, then if $H$ is a hyperplane section of
$E$, the exact sequence
$$\exact{\T_E\otimes\omega_E^{-n+1}}{\T_E\otimes\omega_E^{-n}}
{(\T_E\otimes\omega_E^{-n})|_H}$$ shows us that $H^1(\T_E\otimes\omega^{-n}_E)=0$
for all $n\ge 1$. Now $E$ is $\Ptwo$, $\Pone\times\Pone$, or $\Ptwo$ blown up in
$9-k$ points. In the first two cases, $H^1(\T_E)=0$ is immediate. In the other
cases, $H^1(\T_E)$ is the tangent space to $\Def(E)$, and the moduli space of
non-singular del Pezzo surfaces of degree $\ge 5$ consists of a single point
(reduced since $H^2(\T_E)\cong \dual{H^0(\Omega^1_E\otimes\omega_E)}=0$) by \cite[VII, 2]{[5]}.
Thus $H^1(\T_E)=0$. 
\end{proof}

\begin{say}
\label{(5.5)}
We review here the deformation theory of the singularity $(X,0)$ which is
the
cone over a non-singular del Pezzo surface $E$ of degree $k$, $5\le k\le 9$. The
cases
$6\le k\le 9$ follow from \cite{[2]}.
\begin{itemize}
\item[$k=9$:] $(X,0)$ is rigid. (This also follows from \cite{[42]}.)
\item[$k=8$:] There are two cases. If $E\cong\Pone\times\Pone$, then $(X,0)$ can be
smoothed by taking a hyperplane section of a cone over $\Pthree$ embedded via the
2-uple embedding. If $E\cong F_1$, there is no smoothing. In both cases, $\dim_k
T^1=1$. 
\item[$k=7$:] $(X,0)$ can be smoothed by taking a hyperplane section of a cone over
$\Pthree$ blown up at a point, embedded in $\P^8$ by projecting $\Pthree$ embedded
in $\P^9$ via the 2-uple embedding from a point on the $\Pthree$. Here $\dim_k
T^1=2$.
\item[$k=6$:] There are two distinct smoothings, one coming from taking a hyperplane
section of a cone over $\Pone\times\Pone\times\Pone\subseteq \P^7,$ the other
from taking two hyperplane sections of a cone over $\Ptwo\times\Ptwo\subseteq\P^8$.
Here $\dim_k T^1=3$.
\item[$k=5$:] Any codimension 3 Gorenstein subscheme of the spectrum of a regular
local ring is a Pfaffian subscheme \cite{[3]}, and any Pfaffian subscheme is
smoothable by \cite{[22]}. Here
$\dim_k T^1=4$.
\end{itemize}
\end{say}

We need one more technical lemma:

\begin{lemma} 
\label{5.6} 
Let $E$ be a del Pezzo surface which is either rational
and normal of any degree or else
is non-normal of degree 7 of type (ii) (a) or (b) as given in the proof of Theorem
\ref{5.2}.
Then
\begin{itemize}
\item[(i)] $H^2(\T_E\otimes \omega_E^{-n})=0$ for all $n\ge 0$.
\item[(ii)] $H^2(\Omega^1_E\otimes \omega_E^{-n})=0$ for all $n\ge 0$.
\item[(iii)] $E$ is smoothable.
\end{itemize}
\end{lemma}

\begin{proof}
(i) By Serre duality, $\dual{H^2(\T_E\otimes\omega_E^{-n})}\cong H^0(
\dual{\dual{\Omega^1_E}}\otimes\omega_E^{n+1})$. If $E$ is a normal
del Pezzo surface and thus has only quotient singularities, then
$H^0(\dual{\dual{\Omega^1_E}})=0$ by Hodge theory. Since $\omega_E^{-1}$ is
effective, $H^0(\dual{\dual{\Omega^1_E}}\otimes \omega_E^{n+1})=0$ for all 
$n\ge 0$ also.

If $E$ is non-normal, let $n:\tilde E\rightarrow E$ be the normalization, so 
that $\tilde E$ is a scroll. Then the map $n^*\Omega^1_E\rightarrow
\Omega^1_{\tilde E}$ yields a map $\dual{\dual{\Omega^1_E}}\rightarrow
\dual{\dual{(n_*\Omega^1_{\tilde E})}}\cong n_*\Omega^1_{\tilde E}$,
as is seen by applying duality for a finite morphism twice 
(\cite[III Ex. 6.10 b)]{[16]}
and  \cite{[39]}).
Since
$\dual{\dual{\Omega^1_E}}$ is torsion-free, this map is injective, so
$H^0(\dual{\dual{\Omega^1_E}})\subseteq H^0(\Omega^1_{\tilde E})=0$. Thus the
result still follows in this case.

(ii) If $E$ is normal, then again
$H^2(\Omega^1_E)\cong H^2(\dual{\dual{\Omega^1_E}})=0$ by Hodge theory, and so (ii)
follows. Now suppose $E$ is not normal, with $n:\tilde E\rightarrow E$ the
normalization. We need to consider type (ii) (a) and (b) separately. If $E$ is of
type (a), $n$ maps the section $C_0$ two-to-one to the singular curve $l$. Let
$\iota:C_0\rightarrow C_0$ be the induced involution on $C_0$. $\iota$ 
induces an involution $\iota^*:n_*\Omega^1_{C_0}\rightarrow n_*\Omega^1_{C_0}$.
Let $\shF\subseteq n_*\Omega^1_{C_0}$ be the sheaf of anti-invariants of this
involution. $\shF$ is the image of the map
$\delta':n_*\Omega^1_{C_0}\rightarrow n_*\Omega^1_{C_0}$ given by $\delta'(\alpha)=
\alpha-\iota^*(\alpha)$. Since $n_*\Omega^1_{C_0}=n_*\O_{C_0}(-2)=\O_l(-2)\oplus
\O_l(-1)$, in fact $\shF\cong\O_l(-1)$. Let $\tau^1_E\subseteq\Omega^1_E$ be the
torsion subsheaf of $\Omega^1_E$. We then have a complex
\begin{equation}
\label{(5.7)}
0\mapright{} \tau^1_E \mapright{}\Omega^1_E\mapright{} n_*\Omega^1_{\tilde E}
\mapright{\delta} \shF\mapright{} 0
\end{equation}
where $\delta$ is the (surjective) composition of the surjective restriction map
$n_*\Omega^1_{\tilde E}\rightarrow n_*\Omega^1_{C_0}$ and $\delta'$. This sequence
is exact at $\Omega^1_E$ since $n_*\Omega^1_{\tilde E}$ is torsion-free. It is
exact at $n_*\Omega^1_{\tilde E}$ where $E$ has normal crossings by \cite[(1.5)]{[7]},
and thus this complex splits into exact sequences
$$\exact{K}{n_*\Omega^1_{\tilde E}}{\shF}$$
and 
$$\exact{\Omega^1_E/\tau^1_E}{K}{\tau'}$$
with $\tau'$ supported on points, from which we conclude that $H^2(\Omega^1_E)=0$
as desired, and (ii) follows in this case. 

If $E$ is of case (ii) (b), then we can follow a similar procedure. The
normalization map $n$ maps $C_0$ and $f$ to the singular line $l$. Let $\iota:
C_0\coprod f\rightarrow C_0\coprod f$ be the induced involution interchanging
these two lines. This induces $\iota^*: n_*\Omega^1_{C_0\coprod f}
\rightarrow n_*\Omega^1_{C_0\coprod f}$; let $\shF\cong \Omega^1_l$ be the
anti-invariant part. Again, $\shF$ is the image of $\delta':n_*\Omega^1_{C_0\coprod f}
\rightarrow n_*\Omega^1_{C_0\coprod f}$ given by
$\delta'(\alpha)=\alpha-\iota^*(\alpha)$. We still have the complex 
\eqref{(5.7)}, 
and we finish the argument as before, noting now that the map $H^1(\Omega^1_{\tilde
E})\rightarrow H^1(\shF)$ is surjective.

(iii) A nodal del Pezzo surface is easily seen to be smoothable: the local-global $\ext$
sequence yields
$$\ext^1(\Omega^1_E,\O_E)\rightarrow H^0(\lext^1(\Omega^1_E,\O_E)) \rightarrow
H^2(\T_E)=0,$$ 
and since $\Def(E)$ is smooth, any small deformation of a neighborhood of the
rational double points of $E$ is realised by a global deformation of $E$.

If $E$ is not normal, then $E$ is degree 7 of type (ii) (a) or (b). First suppose
$E\subseteq\P^7$ is of type (ii) (a). Then there is a map $p:E\rightarrow l$
taking a point $x\in E$ to the point of $l\subseteq E$ which is contained in the
same ruling of $E$ as $x$. Each fibre is a singular conic: a union of two $\Pone$'s
or in two cases a doubled line. Each conic spans a plane, and thus $E$ is contained
in a three dimensional scroll. Abstractly, this scroll can be described as the image
of the $\Ptwo$-bundle $\P(p_*\O_E(1))$ in $\P^7$. Note that $p_*\O_E(1)$ is
generated by global sections since $\O_E(1)$ is. Since
$h^0(p_*\O_E(1))=h^0(\O_E(1))=8$, $p_*\O_E(1)$ must be a rank $3$ vector bundle on
$l$ isomorphic to
$\E_{(a,b,c)}=\O_{\Pone}(a)\oplus\O_{\Pone}(b)\oplus\O_{\Pone}(c)$ with
$a,b,c\ge 0$ and $a+b+c+3=8$.
If $t=c_1(\O_{\P(\E_{(a,b,c)})}(1))$ and $f$ is the class of a fibre
of $\P(\E_{(a,b,c)})$, then $\Pic\P(\E_{(a,b,c)})=\ZZ t\oplus\ZZ f$. In order
for $E\subseteq \P(\E_{(a,b,c)})$ to be a conic bundle of degree $7$, it must
have class $2t-3f$. The linear system $|t|$ induces the map of $\P(\E_{(a,b,c)})$
into $\P^7$, and if $|t|$ were not ample (i.e. one of $a$, $b$, or $c$ were zero),
it is easy to see that this would imply two different fibres of $p$ were not
disjoint.
Thus $(a,b,c)=(1,2,2)$ or $(1,1,3)$. In the latter case, the linear system
$|2t-3f|$ has a fixed component given by $t-3f$, and so
$E\subseteq\P(\E_{(1,2,2)})$.

Now a non-singular del Pezzo surface $E'$ of degree 7 also has a conic bundle
structure, so an identical argument also shows that $E'\subseteq\P(\E_{(1,2,2)})$
and is in the same linear system $|2t-3f|$. Thus $E$ is smoothable.

If $E$ is of type (ii) (b), it is easy to see that it is a degenerate case of (ii)
(a),   and so this case is also  smoothable. One can consult \cite{[51]} for an
explicit description of this degeneration. 
\end{proof}

\begin{theorem}
\label{5.8} 
Let $\pi:\tilde X\rightarrow X$ be a primitive type II
contraction with exceptional divisor $E$. Then $X$ is smoothable unless 
\begin{itemize}
\item[(1)] $E\cong\Ptwo$ or
\item[(2)] $E\cong F_1$. 
\end{itemize}
\end{theorem}

\begin{proof}
We want to apply Theorem \ref{4.3}, so we need to verify the hypotheses of the
theorem. First, since $\pi$ corresponds to the contraction of an extremal ray
corresponding to $K_{\tilde X}+E$, $X$ is $\QQ$-factorial by \cite[Prop.~5.1.6]{[20]}.

From now on we will assume that $k\ge 5$. If $k\le 4$, then $X$ has only complete
intersection singularities and we can apply Theorem \ref{3.8}. We have to show that the
singularity of $X$ is good.

Let $P\in X$ be the singular point, with $(\tilde X,E)\rightarrow (X,P)$ the
resolution of the germ $(X,P)$. We first consider deformations of the 
inclusion map $i:E\rightarrow (\tilde X,E)$ of $E$ into
the germ $(\tilde X,E)$. If we denote by $T^1_i$ the tangent
space to the deformation space of the triple $((\tilde X,E), E, i)$, there is an
exact sequence by \cite{[32]}
\begin{align*}
&\hom(\Omega^1_E,\O_E)\oplus H^0(\T_{\tilde X})\mapright{}
H^0(\T_{\tilde X}|_E)\mapright{}T^1_i\\
\mapright{(d_1,d_2)}&
\ext^1(\Omega^1_E,\O_E)\oplus H^1(\T_{\tilde X})
\mapright{} H^1(\T_{\tilde X}|_E).
\end{align*}
Here, by $H^i(\T_{\tilde X})$, we mean cohomology on the germ $(\tilde X,E)$.
Note that $\tilde X$ is non-singular but $E$ may not be. The composed maps 
$d_1:T^1_i\rightarrow\ext^1(\Omega^1_E,\O_E)$ and $d_2:T^1_i\rightarrow
H^1(\T_{\tilde X})$ are the differentials of the maps 
$\Def((\tilde X,E),E,i)\rightarrow \Def(E)$ and $\Def((\tilde X,E),E,i)\rightarrow
\Def(\tilde X,E)$ respectively. 

{\it Claim 1:} $d_1$ and $d_2$ are surjective.

\begin{proof}
Note that $d_2$ is surjective if the map $\ext^1(\Omega^1_E,\O_E)\rightarrow
H^1(\T_{\tilde X}|_E)$ is surjective, which is the case since
$H^1(\dual{(\I_E/\I_E^2)})= H^1(\omega_E)=0$. (In fact, since
$H^0(\dual{(\I_E/\I_E^2)})=0$ also, this map is an isomorphism.) Similarly,
$d_1$ is surjective if
$H^1(\T_{\tilde X})\rightarrow H^1(\T_{\tilde X}|_E)$ is surjective, which is true if
\[
H^2(\T_{\tilde X}(-E))=\lim_{\leftarrow} H^2(\T_{\tilde X}(-E)\otimes\O_{\tilde
X}/\I_E^n)=0.
\]
From
$$\exact{\T_E(-E)}{\T_{\tilde X}(-E)|_E}{\O_E}$$
we see that $H^2(\T_{\tilde X}(-E)|_E)=0$ if
$H^2(\T_E(-E))=H^2(\T_E\otimes\omega_E^{-1})=0$. This is the case by Lemma \ref{5.6} (i).
Tensoring the exact sequence
$$\exact{\I_E^n/\I_E^{n+1}=\omega_E^{-n}}{\O_{\tilde X}/\I_E^{n+1}}{\O_{\tilde X}
/\I_E^n}$$ 
with $\T_{\tilde X}(-E)$ shows that $H^2(\T_{\tilde X}(-E)\otimes\O_{\tilde
X}/\I_E^n)=0$ for all $n$, so $d_1$ is surjective. 
\end{proof}

{\it Claim 2:} $\Def((\tilde X,E), E, i)$ is unobstructed.

\begin{proof}
As noted in the proof of Claim 1, the map $\ext^1(\Omega^1_E,\O_E)
\rightarrow H^1(\T_{\tilde X}|_E)$ is an isomorphism. Also the map
$\hom(\Omega^1_E,\O_E)\rightarrow H^0(\T_{\tilde X}|_E)$ is a surjection. This
shows that $T^1_i\cong H^1(\T_{\tilde X})$, and it then follows from Proposition
\ref{3.4} as in the proof of \cite[Thm.~2.1]{[35]}, that $\Def((\tilde X,E), E, i)$ is
unobstructed. 
\end{proof}

Thus, since $d_2$ is surjective by Claim 1, any
deformation of $(\tilde X,E)$ is of the form $(\tilde X',E')$ with $E'$
a deformation of $E$. Since $d_1$ is surjective, for any small deformation of $E$
to $E'$, there is a deformation $(\tilde X',E')$ of $(\tilde X,E)$. Thus in
particular, if $E$ is smoothable, the general deformation of $(\tilde X,E)$
to $(\tilde X',E')$ yields $E'$ smooth. By Lemma \ref{5.6} (iii),
any of the possible del
Pezzo surfaces under consideration are smoothable.

To summarize, $(\tilde X,E)$ can be deformed to $(\tilde X',E')$ with $E'$
non-singular. So $(X,P)$ can be deformed to $(X',P')$ with $(\tilde
X',E')\rightarrow
(X',P')$ the blow-up, and by Theorem \ref{5.4}, $(X',P')$ is analytically isomorphic to
a cone over a del Pezzo surface. By \ref{(5.5)}, $(X',P')$ is smoothable unless
$E'\cong\Ptwo$ or $F_1$. In this latter case, as it is easy to see that the only
normal del Pezzo surfaces of degree $8$ or $9$ are $\Ptwo$, $F_1$,
$\Pone\times\Pone$ or a quadric cone appropriately embedded, we must also have
$E\cong\Ptwo$ or $F_1$.

If $E$ is not $\Ptwo$ or $F_1$, then $(X',P')$ is smoothable. Thus the general
point in the image of $\Def(\tilde X,E)\rightarrow \Def(X,P)$ is smoothable, so this
image is contained in a smoothing component. The only remaining thing to check is
that $\Def(\tilde X,E)\rightarrow \Def(X,P)$ is an immersion. Thus we need to show
that the differential of this map is injective. This differential
is given by Lemma
\ref{3.2} to be the map
$H^0(R^1\pi_*\T_{\tilde X})\rightarrow H^2_Z(\T_X)$, and the kernel of this map is
$H^1_E(\T_{\tilde X})\cong
\dual{H^0(R^2\pi_*\Omega^1_{\tilde X})}$, which is
easily seen to be zero using similar methods as above via Lemma \ref{5.6} (ii).
\end{proof}

\end{document}